\documentclass[preprint,review,10pt]{elsarticle}

\usepackage[left=3cm, right=3cm, top=3cm, bottom=3cm]{geometry}

\journal{Finite Elements in Analysis and Design}

\usepackage{amsfonts}
\usepackage{bbm}
\usepackage{amssymb}

\usepackage[english]{babel}
\usepackage{parskip}
\usepackage{amsmath}
\usepackage{bm}
\usepackage[intoc]{nomencl}
\usepackage{algorithmic}
\usepackage{algorithm}
\usepackage{nicefrac}
\usepackage{mathtools}
\usepackage{xcolor}
\usepackage[framemethod=tikz]{mdframed}
\usepackage{booktabs}
\usepackage[toc,page]{appendix}
\usepackage{subcaption}
\usepackage{graphicx}
\usepackage{array}
\usepackage{xcolor}
\usepackage{soul}
\usepackage[colorlinks=true,linkcolor=red]{hyperref}
\usepackage{cleveref}
\usepackage{placeins}
\usepackage{xspace}
\usepackage{tikz}
\usepackage{pgfplots}
\usepackage{cancel}
\usepackage{ulem}

\newcommand{\bs}{\boldsymbol}

\crefname{appsec}{Appendix}{Appendices}

\biboptions{numbers,sort&compress}

\begin{document}

\begin{frontmatter}

\title{Elasto-plastic large deformation analysis of multi-patch thin shells by isogeometric approach}

\author[ath1]{G. D. Huynh}
\author[ath1]{X. Zhuang\corref{cor1}}
\author[ath2]{H. G. Bui}
\author[ath2]{G. Meschke}
\author[ath3,ath4]{H. Nguyen-Xuan}

\address[ath1]{Institute for Continuum Mechanics, Leibniz University Hannover, Germany}
\address[ath2]{Institute for Structural Mechanics, Ruhr University Bochum, Germany}
\address[ath3]{CIRTech Institute, Ho Chi Minh City University of Technology (HUTECH), Ho Chi Minh City, Vietnam}
\address[ath4]{Department of Architectural Engineering, Sejong University, 98 Kunja Dong, Kwangjin Ku, Seoul, 143-747, South Korea}
\cortext[cor1]{Corresponding author}

\section*{Highlights}
\begin{itemize}
	\item  A unified thin shell formulation is established, allowing arbitrary nonlinear material models and multi-patch shell structure to be applicable.
	\item The Kirchhoff-Love shell theory is employed and the $C^1$ continuity at patch boundaries is remained by the bending strip method.
	\item The B\'ezier decomposition concept is used to retain the local support of the traditional finite element.
	\item  The computational model is validated and its numerical results agree well with ones in literature.
\end{itemize}

\begin{abstract}
This paper studies elasto-plastic large deformation behavior of thin shell structures using the isogeometric computational approach with the main focus on the efficiency in modelling the multi-patches and arbitrary material formulations. In terms of modelling, we employ the bending strip method to connect the patches in the structure. The incorporation of bending strips allows to eliminate the strict demand of the $C^1$ continuity condition, which is postulated in the Kirchhoff-Love theory for thin shell, and therefore it enables us to use the standard multi-patch structure even with $C^0$ continuity along the patch boundaries. Furthermore, arbitrary nonlinear material models such as hyperelasticity and finite strain plasticity are embedded in the shell formulation, from which a unified thin shell formulation can be achieved. In terms of analysis, the B\'ezier decomposition concept is used to retain the local support of the traditional finite element. The performance of the presented approach is verified through several numerical benchmarks.

\end{abstract}

\begin{keyword}
Isogeometric analysis, Kirchhoff-Love shell theory, Multi-patch structures, Finite Strain, B\'ezier decomposition.
\end{keyword}

\end{frontmatter}

%

\textcolor{blue}{
 \textbf{Nomenclatures}
}

\begin{center}
	\begin{tabular} {>{\color{blue}}l >{\color{blue}}l  }
 \textbf{u}  & \hspace{2cm}  Displacement vector  \\
  \textbf{x}  & \hspace{2cm}  Position Vector in reference configuration  \\
   \textbf{X}  & \hspace{2cm}  Position Vector in current configuration  \\
  $\boldsymbol{\xi}$   & \hspace{2cm} Curvilinear coordinates  \\
  $\boldsymbol{G}_{\alpha}$  & \hspace{2cm} Covariant base vectors in reference configuration \\
  $\boldsymbol{G}^{\alpha} $  & \hspace{2cm} Contravariant base vectors in reference configuration \\
  $G_{\alpha\beta} $  & \hspace{2cm} Covariant metric coefficients \\
  $B_{\alpha\beta} $  & \hspace{2cm} Covariant curvature coefficients \\
  $\boldsymbol{F} $  & \hspace{2cm} Deformation gradient tensor \\
    $\boldsymbol{F}^e $  & \hspace{2cm} Elastic part of the deformation gradient tensor \\
      $\boldsymbol{F}^p $  & \hspace{2cm} Plastic part of the deformation gradient tensor \\
   $\boldsymbol{C} $  & \hspace{2cm} Right Cauchy deformation tensor \\
      $\boldsymbol{C}^p $  & \hspace{2cm} Plastic right Cauchy deformation tensor \\
      $\boldsymbol{b}^e $  & \hspace{2cm} Elastic left Cauchy deformation tensor \\
      $\lambda_{33} $  & \hspace{2cm} thickness stretch \\
     $\boldsymbol{E} $  & \hspace{2cm} Green-Lagrange strain tensor \\
     $\boldsymbol{S} $  & \hspace{2cm} 2nd Piola-Kirchhoff stress tensor \\
     $\varepsilon_{\alpha\beta} $  & \hspace{2cm} Membrane strain components \\
     $\kappa_{\alpha\beta} $  & \hspace{2cm} Bending strain components \\
      $\Psi^{e} $  & \hspace{2cm} Strain energy density \\
     $\Psi^{vol}, \, \Psi^{dev} $  & \hspace{2cm} Volumetric and deviatoric parts of the strain energy density \\
     $\alpha $  & \hspace{2cm} Equivalent plastic strain \\
     $k(\alpha) $  & \hspace{2cm} Hardening function \\
     $\Phi $  & \hspace{2cm} Yield function \\
     $\boldsymbol{\mathbbm{C}} $  & \hspace{2cm} Material tensor \\
          $\boldsymbol{n} $  & \hspace{2cm} Membrane force vector \\
               $\boldsymbol{m} $  & \hspace{2cm} Bending moment vector \\

	\end{tabular}
\end{center}

\section{Introduction}

Thin shell structures play a vital role in the automotive and aerospace industry. Due to the efficiency and the accuracy of the shell formulation, it is typically used to simulate the manufacturing processes or to determine the weak spot of the structure during operation. To accommodate for higher computational fidelity, large strains and rotations are used in the shell formulation. Because the material composing shell is typically metal, the elastoplastic constitutive model is used to investigate the failure mode, e.g. strain localization, under large loading.

There are two typical shell models accounting for the inelastic response of structures. The first model relies completely on stress resultants in which the normal vector of the shell is assumed to be inextensible. The description of the shell surface response is embodied in constitutive equations and the spread of the plasticity through the shell thickness can not be reached completely. Accordingly, formulations of stress resultant constitutive relations are awkward to derive and also sophisticated to implement in the finite element framework, for further discussion see \cite{Dujc2012}. In the standard model, integration points through the thickness direction of the shell are defined in association with the notion of stress-based elastoplasticity and the distribution of plasticity through the thickness can be represented. Therefore, the stress-based approach is more preferred in the finite element context. The model is originally formulated from the three dimensional theory, which can be compatible with formulations of shell kinematics in a direct or indirect manner. Solid-shell element as presented in \cite{Caseiro2015} , where only displacement degree-of-freedom are used, can match perfectly with 3D constitutive model without the need of ad hoc assumptions. However, this element is usually implemented together with methods such as assumed natural strain (ANS) and enhanced assumed strain (EAS) to alleviate locking effects. In contrast, elements based on Kirchhoff-Love hypothesis, e.g., \cite{Kiendl2009} are though rotation-free, need to resort to plane stress enforcement to take through-thickness behaviours into account. The implementation of three dimensional constitutive models for finite strain plasticity with the plane stress assumption can be classified into two common approaches: plane stress-projected constitutive models and the nested iteration approach for the plane stress enforcement, which are also applicable to hyperelasticity, see details in \cite{Neto2008} .The former approach is only applicable if involved equations of the models are adequately simple so that transverse components can be eliminated from the formulation, while the latter one can be implemented in a straightforward manner for any models and is adopted in this work, see also its applicability in recent works \cite{Kiendl2015,Ambati2018}.

Isogeometric analysis (IGA) emerges as a versatile tool to perform analysis and modeling simultaneously using the same interpolation bases.  IGA offers numerical properties which are often beneficial in numerical analysis, e.g. positiveness of the basis functions, higher continuity. IGA is shown to be superior over standard FEM based on Lagrange shape function for dynamics analysis . Refinement with IGA is more straightforward than with the standard FEM counterpart, in which hpk-refinement scheme can be used to offer both order elevation and addition of degree of freedoms (d.o.fs) conveniently without introducing internal modes or hanging nodes. It is noted that, refinement in IGA preserves the geometry, and thus eliminating the geometrical error stemming from modelling. On the other hand, IGA is particularly suited for shell modelling. The smoothness of NURBS basis functions and the similarity between natural coordinates of the isogeometric parameterization and the shell coordinates lead the displacement based shell formulations to be implemented in a direct manner without the need for coordinate transformation. In addition, NURBS is also utilized as surface modelling tool in the CAD technology. This suggests that, extension of shell formulation for IGA is practical, and is apparently beneficial for many industrial applications. Nevertheless, many complex shell structures are not easy to be modelled with a single NURBS patch. An option is to use a modelling technique that supports connecting control point, i.e. T-splines . However, this technique is not broadly implemented in many CAD packages. It is therefore the preferred option to use multi NURBS patch to model a shell structure. \textcolor{blue}{In terms of shell formulations, the Mindlin-Reissner theory is commonly-used in the finite element context as only C0-continuity is required. Typical shell elements for this theory are generations of Mixed Interpolation of Tensorial Components (MITC) proposed by K. Bathe et al. \cite{Bathe1984,Bathe2016, Bathe2017}, shell elements based on discrete Kirchhoff-Love constraints by Areias et al. \cite{Areias2005a} \cite{Areias2005b}. These works need adhoc techniques such as an enhanced assumed strain to alleviate locking effects that are the most challenged hinge in the application of the finite element method (FEM) to shell problems. For a comprehensive overview of other FEM shell elements, readers are referred to an excellent work .In contrast to the Mindlin-Ressner theory,} shell formulations using the Kirchhoff-Love hypothesis require $C^1$ continuity of the basis functions to ensure the smoothness of bending moment interpolation. It is well known that, if the open knot vector is used in the modelling, only $C^0$ continuity is achieved along the patch boundary. Thus, additional technique, i.e. bending strip method, mortar method, shall be used to maintain the $C^1$ patch continuity. While the mortar method or Nitsche method that requires the constraints on the derivatives of the basis function at the boundary, the bending strip method only adds fictitious bending stiffness on the overlapping domain in the vicinity of the boundary and its term does not complicate the weak form and can be evaluated locally. Therefore, the latter is simple to implement, yet ensures the reliability in the description of structure behaviours.

In this paper, an efficient computational model for thin shell using isogeometric analysis is presented. In terms of modelling, the efficiency is preserved by using multi NURBS patch structure with $C^0$ continuity on the patch boundaries. The bending strip method \citep{Kiendl2010b} is employed to maintain the $C^1$ continuity condition of the Kirchhoff-Love formulation. Moreover, arbitrary 3D nonlinear materials such as hyperelasticity, finite strain plasticity can be included in the formulation by enforcing the plane stress condition. This leads to our main focus in this work in which a unified thin shell formulation is established to  not only model shell structures with complex geometries or multiple patches, but also exhibit a wide range of material behaviours. Note that, for the plane stress enforcement, we follow the idea in \cite{Neto2008}, in which the plane stress constraint is considered an additional equation that is used to determine the in-plane stress components and the thickness stretch. Though similar approaches are presented in recent works \cite{Kiendl2015} \cite{Ambati2018}, these methods are restricted to single-patch shell geometries and can not be applicable to a large class of shell structures involving multi-patches.
On the aspect of analysis, the B\'ezier decomposition method \citep{Borden2011} is employed. The advantage of the B\'ezier decomposition approach is twofold: 1) maintaining the local characteristic of the finite element and improving the performance of the global stiffness matrix assembly procedure, 2) enabling the reuse of standard finite element components, including the kinematics and the constitutive law. To verify the applicability of the proposed computational model, it is validated with reference examples using hyperelasticity and elastoplasticity constitutive laws.

 The structure of the paper is as following: in \Cref{sec:theory}, geometry description of the thin shell using Kirchhoff-Love hypothesis is presented, followed by the shell kinematics and the treatment of plane stress enforcement which enables the usage of three dimensional constitutive model for the shell element. \Cref{sec:isogeometric} underlines some important properties of NURBS basis functions and give the introduction to the bending strip method with an illustrative example. In \Cref{sec:plasticity}, required kinematics modification to accommodate for elastoplastic constitutive model in finite strain, which relies on  the hypothesis of multiplicative plasticity framework, is briefly presented. \Cref{sec:discretization} presents the discretization of the weak form based on the directional derivative, and the formation of the B\'ezier finite element, based on B\'ezier decomposition. \Cref{sec:examples} is devoted for selected numerical examples. The paper is concluded by \Cref{sec:conclusion}.


\section{Kirchhoff-Love theory for thin shell}
\label{sec:theory}

\subsection{Notation}

In this paper, the following notation is employed: italic letters ($g$, $G$) are used to indicate scalar and bold letters ($\mathbf{u}$, $\mathbf{g}$, $\mathbf{G}$) indicate the vector and second order tensor. Lower case letters indicate terms in the reference configuration and upper case letters are for terms in the current configuration. $\xi^\alpha$ ($\alpha=1,2$) denote the convective curvilinear coordinates of the shell and $\xi^3$ denotes the local coordinates in thickness direction. The indices in Greek letters ($\alpha$, $\beta$) take values of $\{1, 2\}$ and the indices in Latin letters ($i$, $j$) take values of $\{1, 2, 3\}$.

\subsection{Geometry definition}

The deformation and the rotation of the shell are defined on its mid-surface as
\begin{equation}
\mathbf{u} (\boldsymbol{\xi}) = \mathbf{x} (\boldsymbol{\xi}) - \mathbf{X} (\boldsymbol{\xi})
\label{eq:deformation_rotation}
\end{equation}

where $\mathbf{x}$ and $\mathbf{X}$ are the position vectors of a material point on the mid-surface in the deformed configuration and undeformed configuration. $\boldsymbol{\xi} = \{ \xi^1, \xi^2, \xi^3 \}$ is the local coordinates of the material point.

The covariant tangential base vectors are defined by
\begin{equation}
\mathbf{G}_{\alpha} = \mathbf{X}_{,\alpha} = \dfrac{\partial \mathbf{X}}{\partial \boldsymbol{\xi}^{\alpha}} \qquad \mathbf{g}_{\alpha} = \mathbf{x}_{,\alpha} = \dfrac{\partial \mathbf{x}}{\partial \boldsymbol{\xi}^{\alpha}}
\label{eq:CovariantBases}
\end{equation}

\begin{figure}[htb]
	\centering
	\includegraphics[scale=0.55]{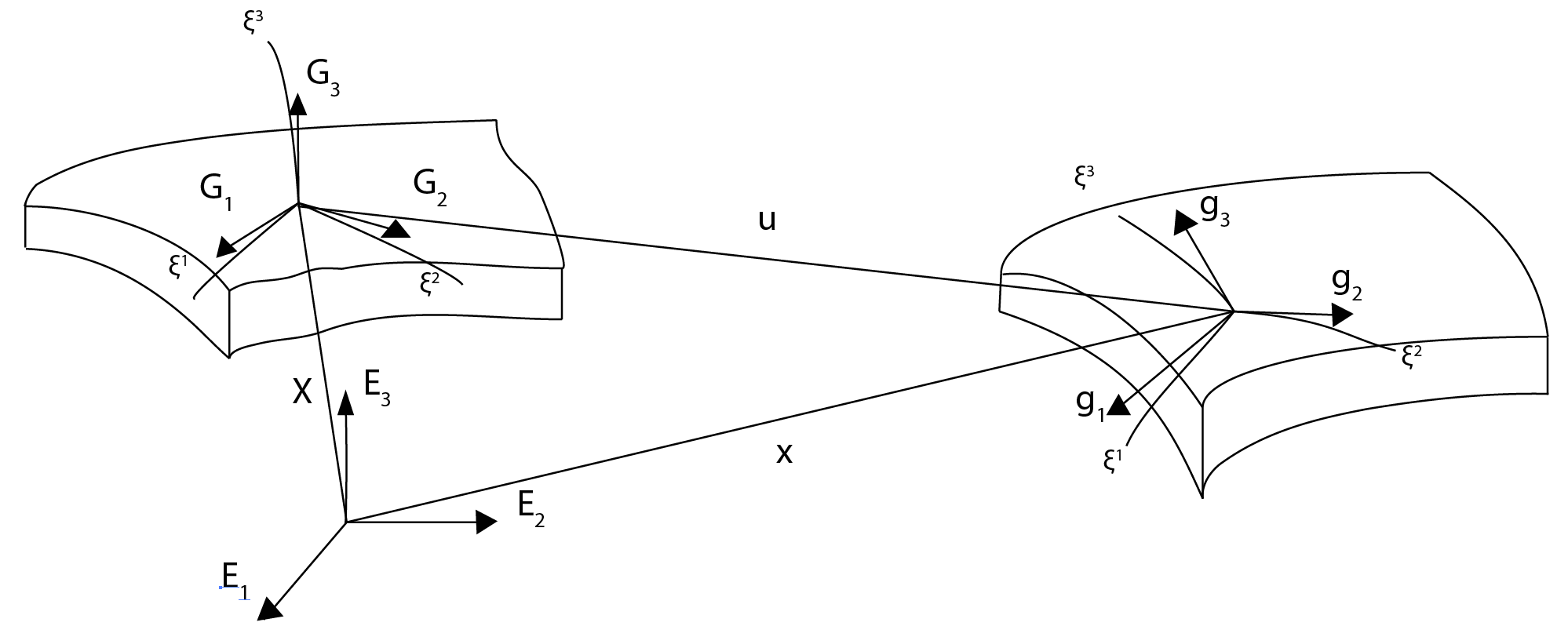}
	\caption{\textcolor{blue}{Description of body motion between different configurations using curvilinear coordinate systems}}
	\label{CurvilinearCoord}
\end{figure}

\textcolor{blue}{Figure. \ref{CurvilinearCoord} reveals the representation of the basis vectors between reference and current configurations in a curvilinear coordinate system.}

The covariant metric coefficients of the surface are computed by
\begin{equation}
G_{\alpha\beta} = \mathbf{G}_{\alpha} \cdot \mathbf{G}_{\beta} \qquad g_{\alpha\beta} = \mathbf{g}_{\alpha} \cdot \mathbf{g}_{\beta}
\label{eq:CovariantMetric}
\end{equation}
and the contravariant base vectors by
\begin{equation}
\mathbf{G}^{\alpha} = G^{\alpha\beta} \mathbf{G}_{\beta} \qquad \mathbf{g}^{\alpha} = g^{\alpha\beta} \mathbf{g}_{\beta}
\label{eq:ContracovariantMetric}
\end{equation}
with ${G^{\alpha\beta}} = ({G_{\alpha\beta}})^{-1}$ and ${g^{\alpha\beta}} = ({g_{\alpha\beta}})^{-1}$.

The unit normal vectors on the middle surface can be computed from the covariant tangential base vectors as
\begin{equation}
\mathbf{G}_3 = \dfrac{\mathbf{G}_1 \times \mathbf{G}_2}{\| \mathbf{G}_1 \times \mathbf{G}_2 \|} \qquad \mathbf{g}_3 = \dfrac{\mathbf{g}_1 \times \mathbf{g}_2}{\| \mathbf{g}_1 \times \mathbf{g}_2 \|}
\label{eq:NormalBases}
\end{equation}

From \Cref{eq:NormalBases,eq:CovariantBases}, one can compute the curvature tensor $\mathbf{B}$ and $\mathbf{b}$. They are defined by the second fundamental form of surfaces
\begin{equation}
\begin{split}
B_{\alpha \beta} &= \dfrac{1}{2} \left( \mathbf{G}_\alpha \cdot \mathbf{G}_{3,\beta} + \mathbf{G}_\beta \cdot \mathbf{G}_{3,\alpha} \right) \\
b_{\alpha \beta} &= \dfrac{1}{2} \left( \mathbf{g}_\alpha \cdot \mathbf{g}_{3,\beta} + \mathbf{g}_\beta \cdot \mathbf{g}_{3,\alpha} \right)
\end{split}
\end{equation}

\subsection{Kinematics}
\label{sec:kinematics}

Strain measurement is represented with respect to the right Cauchy-Green deformation tensor $\bs{C} = \bs{F}^T \bs{F}$, where $\bs{F}$ is the deformation gradient
\begin{equation}
\bs{F} = \frac{\text{d} \mathbf{x}}{\text{d} \mathbf{X}} = \mathbf{g}_i \otimes \mathbf{G}^i.
\label{eq:DeformationGradient1}
\end{equation}

Then, the deformation tensor can be written in terms of contravariant base vectors as
\begin{equation}
\bs{C} = \bs{F}^T \bs{F} = \bs{g}_i \cdot \bs{g}_j  \bs{G}^i \otimes \bs{G}^j  = g_{ij}  \bs{G}^i \otimes \bs{G}^j.
\label{eq:RightDeformationTensor}
\end{equation}

In \Cref{eq:RightDeformationTensor}, components of the right Cauchy-Green deformation tensor are the same as the metric coefficients in the deformed configuration. However, $C_{33}$ that describes the change of the shell thickness during deformations can not take the same value with $g_{33}$ that is equal $1$. Hence, $C_{33}$ will be determined by an additional constraint, not directly by the solution of the midsurface variables. The thickness stretch $\lambda_3$ is computed from $C_{33}$ as
\begin{equation}
\label{ThicknessStretch}
\lambda_3 = \sqrt{C_{33}}.
\end{equation}


The Green-Lagrange strain tensor $\bs{E} = \frac{1}{2} (\bs{C} - \bs{I} )$ is typically separated into a constant part representing membrane action and a linear part representing bending action:
\begin{align}
\bs{E} = E_{ij} \mathbf{G}^i \otimes \mathbf{G}^j \qquad E_{\alpha \beta} = \epsilon_{\alpha \beta} + \xi^3 \kappa_{\alpha \beta}
\end{align}

in which
\begin{align}
\varepsilon_{\alpha \beta} &= \frac{1}{2} (g_{\alpha \beta} - G_{\alpha \beta} ) \\
\kappa_{\alpha \beta} &= b_{\alpha \beta} - B_{\alpha \beta}
\end{align}

In order to be consistent with material models which are often expressed in the Cartesian coordinate frame, the strain quantities and the deformation gradient need to be represented in a local Cartesian coordinate system. A general expression for such a transformation is given by
\begin{equation}
\bar{(\cdot)}_{\gamma \delta} = (\cdot)_{\gamma \delta} (\bs{E}_{\gamma} \cdot \bs{G}^{\alpha}) (\bs{G}^{\beta} \cdot \bs{E}_{\delta}),
\end{equation}
where $\bar{(\cdot)}_{\gamma \delta}$ and $(\cdot)_{\gamma \delta}$ are components of generic second order tensors obtained in local Cartesian and curvilinear coordinates respectively. The local Cartesian basis vectors $\bs{E}_{\gamma}$ and $\bs{E}_{\delta}$ are given by equations
\begin{equation}
\bs{E}_1 = \frac{\bs{G}_1}{||\bs{G}_1||}, \,\,\,\,\,\, \bs{E}_2 = \frac{\bs{G}_2 - (\bs{G}_2 \cdot \bs{E}_1)\bs{E}_1}{||\bs{G}_2 - (\bs{G}_2 \cdot \bs{E}_1)\bs{E}_1||}, \,\,\,\,\,\,   \bs{E}_3= \bs{G}_3
\end{equation}

\subsection{Finite strain kinematics for elastoplasticity}
\label{sec:plasticity}

The kinematics presented in \Cref{sec:kinematics} can be applied for the constitutive laws exhibiting elastic behaviours, such as linear elasticity or hyperelasticity. For constitutive laws accounting for plasticity effects, the elastoplastic kinematics model for finite strain \citep{Simo1988a,Simo1988b} shall be employed, which is based on the concept of the plastically deforming intermediate configuration and the multiplicative decomposition of the deformation gradient.

To account for the plasticity behaviour, the deformation gradient is decomposed as
\begin{equation}
\bs{F} = \bs{F}^e \bs{F}^p
\label{eq:SplitDeformationGradient}
\end{equation}
with $\bs{F}^e$ and $\bs{F}^p$ as the elastic and plastic parts of the deformation gradient tensor.

Following \Cref{eq:SplitDeformationGradient}, the plastic right Cauchy-Green deformation tensor is computed as
\begin{equation}
\bs{C}^p = \left( \bs{F}^{p} \right)^T \bs{F}^p
\label{eq:PlasticRightCauchyGreen}
\end{equation}

and the elastic left Cauchy-Green deformation tensor:
\begin{equation}
\bs{b}^e = \bs{F}^{e} \left( \bs{F}^{e} \right)^T
\label{eq:ElasticLeftCauchyGreen}
\end{equation}

From \Cref{eq:SplitDeformationGradient,eq:PlasticRightCauchyGreen,eq:ElasticLeftCauchyGreen}, the elastic left deformation tensor $\mathbf{b}^e$ can be computed from $\mathbf{C}^p$ and $\mathbf{F}$:
\begin{equation}
\bs{b}^e = \bs{F} \left( \bs{C}^{p} \right)^{-1} \bs{F}^T
\end{equation}

The total rate of the elastic deformation tensor is obtained by Lie derivative as follows
\begin{equation}
\mathcal{L}[\bs{b}^e] = \bs{F} \left( \dot{\bs{C} }^{p} \right)^{-1} \bs{F}^T
\end{equation}

The elastic behaviour of the material is modeled by strain energy density function of hyperelastic material, which is further separated into the deviatoric and volumetric components as
\begin{equation}
\Psi^e = \Psi^e_{vol} + \Psi^e_{dev}
\end{equation}

with
\begin{align}
\Psi^e_{vol} &= \frac{1}{2} K \left ( \frac{1}{2}(\left( J^{e} \right)^2 -1) - \text{ln} J^e \right ) \\
\Psi^e_{dev} &= \frac{1}{2}G\left (  \text{tr}[\hat{\bs{b}}^e] - 3 \right )
\end{align}

where $J^e = \text{det}[\bs{F}^e]$, $\hat{\bs{b}}^e = \left( J^{e} \right)^{-2/3} \bs{b}^e$, and $G$ and $K$ are identified as the shear modulus and the bulk modulus in infinitesimal deformation respectively.

The classical Mises-Huber yield condition established with respect to deviatoric part of the Kirchhoff stress $\bs{\tau}_{dev}$ and the hardening function $k(\alpha)$  is given by
\begin{equation}
\Phi(\bs{\tau}_{dev}, \alpha) = ||\bs{\tau}_{dev}|| - \sqrt{\frac{2}{3}}k(\alpha) \leq 0
\end{equation}
where $\Phi$ is the yield function and $\alpha$ the equivalent plastic strain.
Based on the notion of the principle of maximum plastic dissipation, the associative flow rule  can be expressed by
\begin{equation}
\left( \dot{\bs{C}}^{p} \right)^{-1} = - \frac{2}{3} \gamma \, \text{tr}[\bs{b}^e] \, \bs{F}^{-1} \bs{n} \bs{F}^{-T}
\label{eq:FlowRule}
\end{equation}

with $\bs{n} = \bs{\tau}_{dev}/||\bs{\tau}_{dev}||$ and $\gamma$ is the plastic multiplier. The evolution of the yield stress known as hardening phenomenon is driven by the equation
\begin{equation}
\dot{\alpha} = \sqrt{\frac{2}{3}} \gamma
\label{eq:HardendingEvolution}
\end{equation}

Finally, the so-called loading/unloading conditions of the elastoplastic model are governed by the Karush-Kuhn-Tucker condition:
\begin{equation}
\gamma \geq 0 , \, \, \, \, {\Phi}(\bs{\tau}_{dev}, \alpha) \leq 0, \, \, \, \gamma {\Phi}(\bs{\tau}_{dev}, \alpha) = 0
\end{equation}

accompanied by the consistency condition:
\begin{equation}
\gamma \dot{{\Phi}}(\bs{\tau}_{dev}, \alpha) = 0.
\end{equation}

\subsubsection{Time integration for the finite strain plasticity model}

The solution of \Cref{eq:FlowRule,eq:HardendingEvolution} requires a time integration scheme, for which the backward Euler scheme is adopted, such that
\begin{align}
\left( \bs{C}_{n+1}^{p} \right)^{-1} &= \left( \bs{C}_{n}^{p} \right)^{-1} - \frac{2}{3} \Delta \gamma [\bs{b}^e_{n+1}] \bs{F}_{n+1}^{-1} \bs{n}_{n+1} \bs{F}_{n+1}^{-T} \\
\alpha_{n+1} &= \alpha_n + \sqrt{\frac{2}{3}} \Delta \gamma
\label{eq:EvolutionEquation}
\end{align}
$\Delta \gamma$ is the incremental plastic multiplier and satisfies $\Delta \gamma > 0$.

The classical return mapping scheme \cite{Simo1998} \cite{Simo1988a}  is utilized to integrate \Cref{eq:EvolutionEquation}, and summarized in Table.\ref{tab:ReturnMapping}.
\begin{table}
  \begin{center}
    \caption{Hyperelastic extension of return-mapping algorithm .}
    \label{tab:ReturnMapping}
    \begin{tabular}{l} 
       \toprule[2 pt]
 $\bs{F}_{n+1}$ are given at $t_{n+1}$, $\bs{C}_{n}^{p^{-1}}$ and  $\alpha_n$ are known at $t_{n}$.\\
       Elastic predictor step stemming from intermediate configuration: \\

       \hspace*{5mm} $\bar{\bs{F}}_{n+1} = J_{n+1}^{-1/3} \, \bs{F}_{n+1} $ \\
       \hspace*{5mm} $J_{n+1} = \text{det}[{\bs{F}}_{n+1}] $ \\
       \hspace*{5mm} $\bar{\bs{b}}_{n+1}^{e,trial} = \bar{\bs{F}}_{n+1} \, \bs{C}_n^{p-1} \,   \bar{\bs{F}}_{n+1}^T$ \\
       \hspace*{5mm} $ \bs{\tau}_{dev,n+1}^{trial} = \text{G} \, \text{dev}[\bar{\bs{b}}_{n+1}^{e,trial}] $ \\
       \hspace*{5mm} $ \bs{\tau}_{vol,n+1}^{trial} = \frac{\text{1}}{2}K \, (J^2_{n+1} - 1) \, \bs{I}$ \\
      \hline
      Check for yielding:\\

      \hspace*{5mm} IF $ \Phi_{n+1}^{trial} = ||\bs{\tau}_{dev,n+1}^{trial}|| - \sqrt{\frac{2}{3}} R(\alpha_n) \leq 0$ \\
      \hspace*{8mm} THEN set $(\cdot)_{n+1}  = (\cdot)_{n+1}^{trial}$ and EXIT \\
      \hspace*{5mm} ELSE go to corrector step\\
      \hline
      Return mapping corrector:  \\
     \hspace*{5mm}  Compute $\Delta \gamma $ by solving the equation \\
	\hspace*{5mm} $\Phi_{n+1}^{trial} = ||\bs{\tau}_{dev,n+1}^{trial} || -  \sqrt{\frac{2}{3}} k(\alpha_n +  \sqrt{\frac{2}{3}} \Delta \gamma  ) - \frac{2}{3} \Delta \gamma  \text{G} \text{tr}[\bar{\bs{b}}_{n+1}^{e,trial}] = 0  $ \\
   	\hspace*{5mm} Update the state variables \\

\hspace*{5mm}  $ \bs{\tau}_{dev,n+1} = \bs{\tau}_{dev,n+1}^{trial} - \frac{2}{3} \Delta \gamma  \text{G} \text{tr}[\bar{\bs{b}}_{n+1}^{e,trial}] \bs{n}_{n+1}^{trial}$  \\
\hspace*{5mm} $ \bs{\tau}_{vol,n+1} = \bs{\tau}_{vol,n+1}^{trial}$\\
\hspace*{5mm} $\bs{\sigma}_{n+1} = \frac{1}{J} (\bs{\tau}_{vol,n+1} + \bs{\tau}_{dev,n+1})$ \\
\hspace*{5mm} $\bs{C}_{n+1}^{p^{-1}} = \bs{C}_{n}^{p^{-1}} - \frac{2}{3} \Delta \gamma  \text{tr}[\bar{\bs{b}}_{n+1}^{e,trial} ] \bs{F}_{n+1}^{-1} \bs{n}_{n+1}^{trial} \bs{F}_{n+1}^{-T} $ \\
\hspace*{5mm} $ \alpha_{n+1} = \alpha_n + \sqrt{\frac{2}{3}} \Delta \gamma $\\
     \bottomrule[2 pt]
    \end{tabular}
  \end{center}
\end{table}

\subsection{Zero transverse normal stress enforcement for thin shell}
Here, the indigenous three-dimensional constitutive models are utilized within a Newton-Raphson loop by solving the equation of plane stress constraint at each Gauss point with thickness stretch $C_{33}$ as the unknown.

Following \cite{Neto2008} \cite{Dodds1987}, $C_{33}$ is obtained by satisfying the zero transverse normal stress condition, which takes the form of
\begin{equation}
S^{33} = \mathbbm{C}^{33 \alpha \beta} E_{\alpha \beta} + \mathbbm{C}^{3333} E_{33} = 0
\label{eq:}
\end{equation}

from which the coefficients of the material tensor are modified as
\begin{equation}
\hat{\mathbbm{C}} = \mathbbm{C} - \frac{\mathbbm{C}^{\alpha \beta 33} \mathbbm{C}^{33 \gamma \delta} }{\mathbbm{C}^{3333}}
\label{CondensedTensor}
\end{equation}

$C_{33}$ is solved iteratively by the equation
\begin{equation}
S^{33} + \frac{\partial S^{33}}{\partial C_{33}} \Delta C_{33} = S^{33} + \frac{1}{2} \mathbbm{C}^{3333} \Delta C_{33} = 0
\label{eq:NRPlaneStress}
\end{equation}

followed by  an incremental update
\begin{equation}
\Delta C_{33}^{(I)} = -2 \frac{S^{33,(I)}}{C^{3333,(I)}}
\end{equation}

where $I$ denotes the number of iteration step.

The thickness change is corrected by computing  $C_{33}$ and $\lambda_3$ as follows
\begin{align}
C_{33}^{(I+1)} &= C_{33}^{(I)} + \Delta C_{33}^{(I)} \\
\lambda_3^{(I+1)} &= \sqrt{C_{33}^{(I+1)} }
\end{align}

$\lambda_3$ will be updated until $S^{33}$ converges to a defined tolerance which fulfills the zero transverse normal stress condition. Subsequently, the statically material tensor
$\hat{\mathbbm{C}}$ is obtained and only its in-plane components are chosen for the shell structures.

\subsection{Variational form of the equilibrium equation}

The mechanical response of a structure is governed by the principle of virtual work which requires the sum of internal and external work vanishes in the equilibrium state
\begin{equation}
\label{eq:PVW}
\delta W = \delta W^{int} - \delta W^{ext}
\end{equation}

The internal virtual work is defined by
\begin{align}
\delta W^{int} &=   \int\limits_{\Omega} S_{ij} : \delta E_{ij} d \, \Omega =
 \int\limits_{A} \Big ( n_{\alpha \beta} : \delta  \varepsilon_{\alpha \beta} +  m_{\alpha \beta} : \delta  \kappa_{\alpha \beta} \Big ) \, dA  ,
\end{align}

and the external virtual work takes the form of
\begin{equation}
\delta W^{ext} = \int\limits_{A} \delta u_i f_i \, d A.
\end{equation}

\section{NURBS-based isogeometric analysis for shell structures}
\label{sec:isogeometric}

In this section, the useful properties of NURBS basis functions which are suitable for the analysis of Kirchhoff-Love shells are summarized, followed by an introduction to the bending strip method which is employed in the multi-patch case.

\subsection{Interpolation using NURBS basis functions}
Due to the appearance of  the second derivatives in curvature changes of the shell kinematics, $C^1$-continuity of the employed basis functions must be required to attain a conforming discretization of the Kirchhoff-Love hypothesis. This can be simply achieved using NURBS basis functions which are used to approximate the solution fields and represent the shell geometry exactly. Following that, given a control net $\bs{P}_{i,j},1 \leq i\leq n,1 \leq j \leq m$, polynomial orders $p$ and $q$, and knot vectors $\Xi_1 =\{\xi_1^1, \xi_1^2, ..., \xi_1^{n+p+1}\}$ and $\Xi_2 =\{\xi_2^1, \xi_2^2, ..., \xi_2^{m+q+1}\}$, a shell surface can be represented by sum:
\begin{equation}
\label{SurfaceEquation}
\bs{S}(\xi_1, \xi_2) = \sum \limits_{i=1}^{n} \sum \limits_{j=1}^{m} N_{i,j}(\xi_1, \xi_2) \bs{P}_{i,j}
\end{equation}
where $\xi_1 \in \Xi_1 $, $\xi_2 \in \Xi_2$ are coordinates in the parametric space, $N_{i,j}$ are the NURBS basis functions, which are defined as
\begin{equation}
\label{NURBSBasisFunctions}
N_{i,j}(\xi_1, \xi_2) = \frac{ B_{i,p}(\xi_1)  B_{j,q}(\xi_2) w_{i,j}} { \sum \limits_{k=1}^{n} \sum \limits_{l=1}^{m} B_{k,p}(\xi_1) B_{l,q}(\xi_2) w_{k,l}}.
\end{equation}

In \Cref{NURBSBasisFunctions}, $B$ is the univariate B-spline basis functions and $w$ is the weight factor. In the parametric space, the elements are defined by partitioning the knots into knot spans. Inside each knot span, the NURBS basis functions are $C^{\infty}$-continuous and $C^{p-k}$ at each knot ($k$ denotes the multiplicity of the knot). When $k \leq p-1$, the $C^1$ continuity is satisfied and is useful for the analysis of Kirchhoff-Love thin shell. Nevertheless, in the typical scenario, the open knot vector is used \citep{Hughes2005}, which leads to only $C^0$ continuity at the patch boundaries. Without additional treatment, a kink in bending moment will occur on the patch boundaries if the Kirchhoff-Love formulation  is used. This issue will be further described in the next section.

Displacements of the shell are approximated with respect to the NURBS basis functions which are used to model the midsurface of the shell as
\begin{equation}
\bs{u} = N_I(\xi_1, \xi_2) \bs{u}_I
\end{equation}
where $\bs{u}_I$ are the displacement at the I-th control point. Following that, the first and second derivatives of the displacements can be obtained:
\begin{align}
\bs{u}_{,\alpha} &= \dfrac{\partial N_I(\xi_1, \xi_2)}{\partial \xi_\alpha} \bs{u}_I \\
\bs{u}_{,\alpha\beta} &= \dfrac{\partial^2 N_I(\xi_1, \xi_2)}{\partial \xi_\alpha \partial \xi_\beta} \bs{u}_I
\end{align}

\subsection{The bending strip method}
\begin{figure}[htb]
\centering
\includegraphics[scale=0.3]{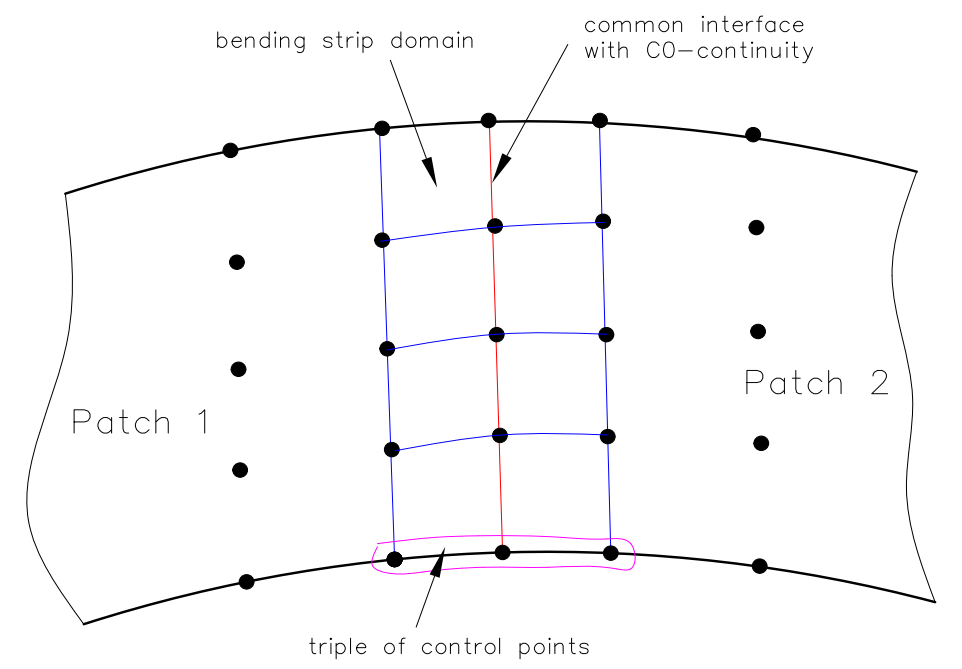}
\quad
\includegraphics[scale=0.4]{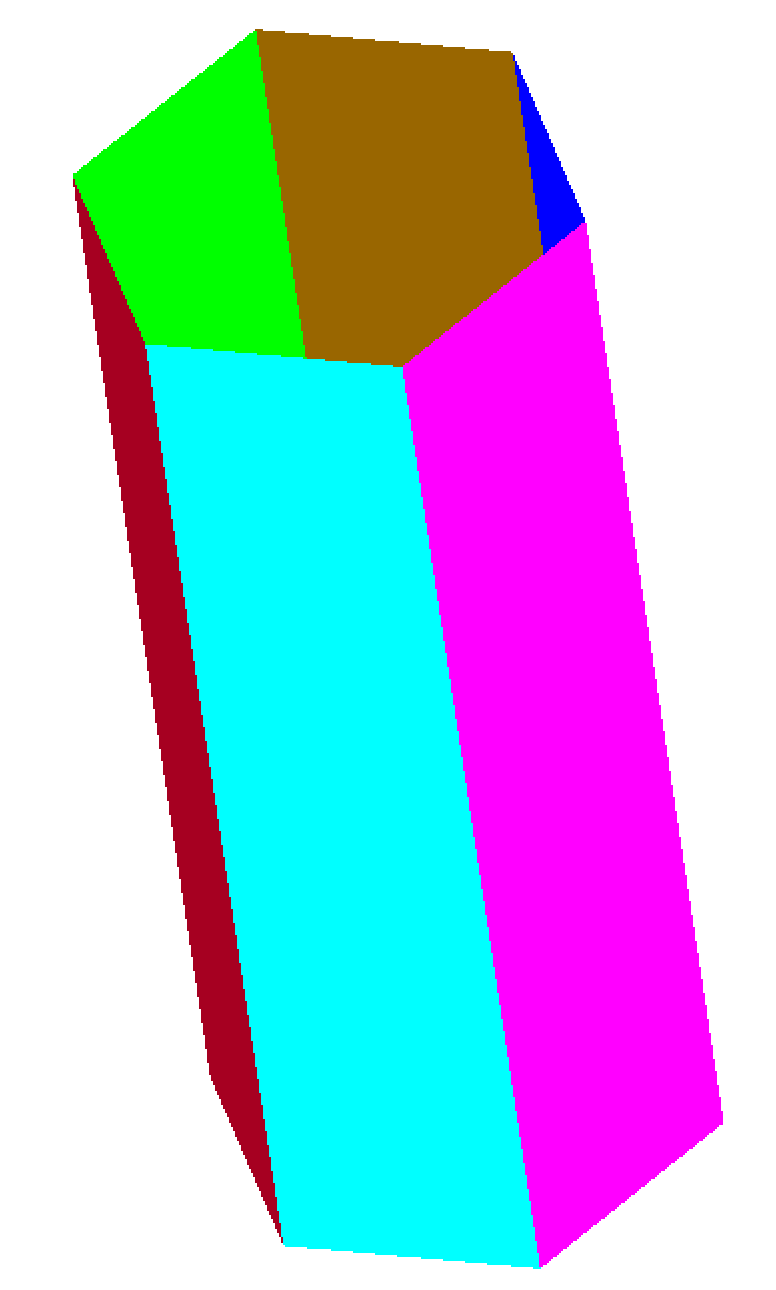}
\caption{Left: Illustration of the bending strip domain; Right: Geometry of the Hex-can example}
\label{StripSchematic}
\end{figure}

The bending strip method \citep{Kiendl2010b} is a mesh connecting method, which fundamentally adds additional kinematics constraints at patch interface in order that the bending moment can be transmitted between patches. The strip is constructed by the triples of control points, one at the shared boundary, one on each side of two involved patches. With this construction, the parametric domain of the bending strip surface comprises of one quadratic element in the direction perpendicular to the interface and linear elements along the boundary. For an illustration of the bending strip method, the reader is referred to \Cref{StripSchematic} (Left). The principle of virtual work in \Cref{eq:PVW} is rewritten to account for the bending strips as
\begin{equation}
\delta W = \delta W^{int} + \delta W^{int}_{strip} - \delta W^{ext}
\label{eq:PVW2}
\end{equation}

The new term $\delta W^{int}_{strip}$ is computed as
\begin{equation}
\delta W^{int}_{strip} = \int\limits_{A} \frac{t^3}{12} \bs{\kappa}^T \bs{\mathbbm{C}}_{strip} \bs{\kappa} \, d A.
\end{equation}

Adopting the bending material stiffness proposed in \cite{Kiendl2010b}, we have
\begin{align}
\bs{\mathbbm{C}}_{strip} =
\begin{bmatrix}
0 & 0 & 0 \\
0 & E_{strip}  & 0  \\
0 & 0 & 0
\end{bmatrix}
\end{align}

in which  the value of the bending strip modulus $E_{strip} $ must be chosen to be high enough and a range of its values varies from $10^3 \times E$ to $10^5 \times E$  with $E$ as Young's modulus.


\begin{figure}[htb]
\centering
  \begin{tabular}{@{}cccc@{}}
    \includegraphics[width=0.4\textwidth]{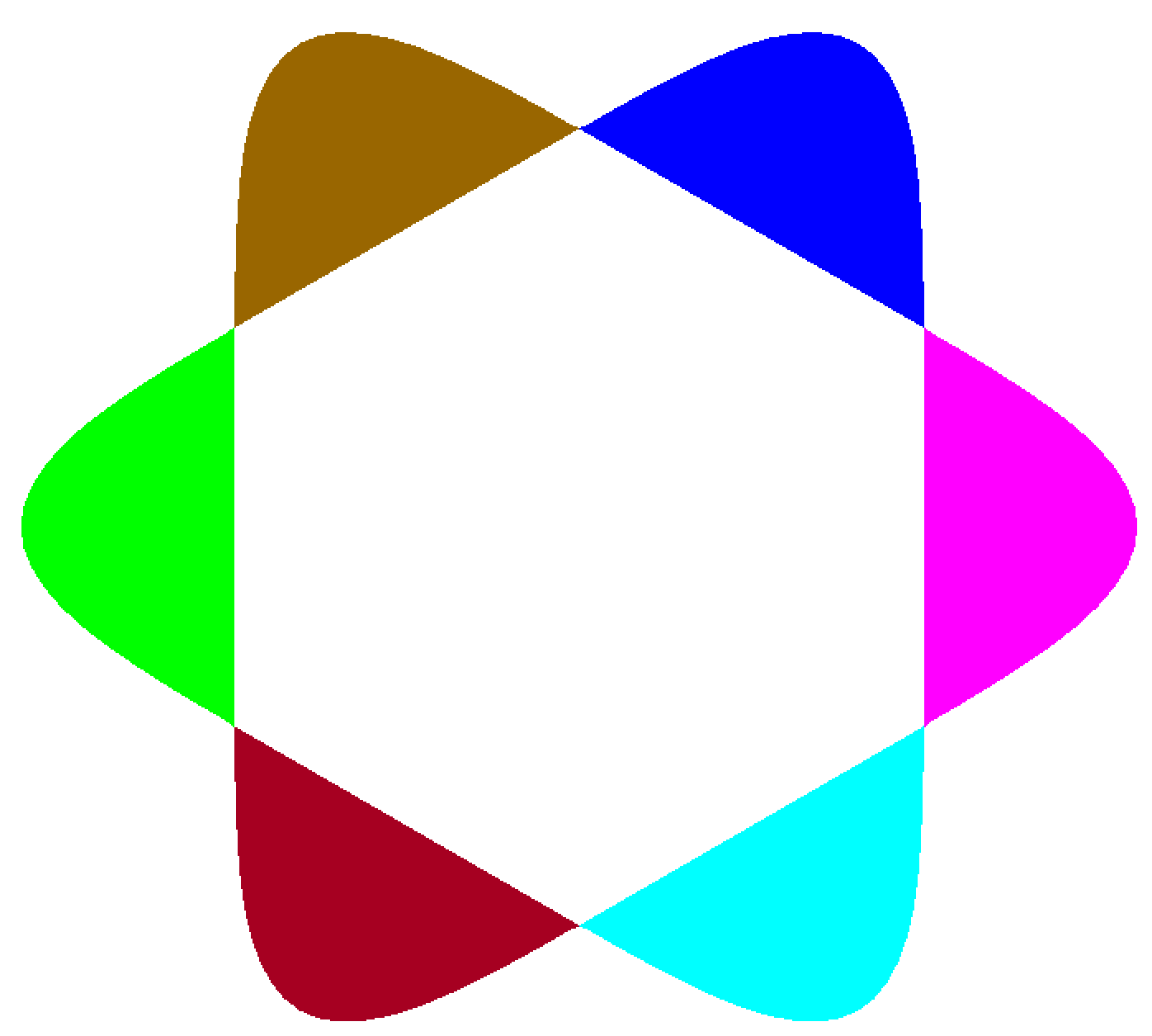}
    \includegraphics[width=0.4\textwidth]{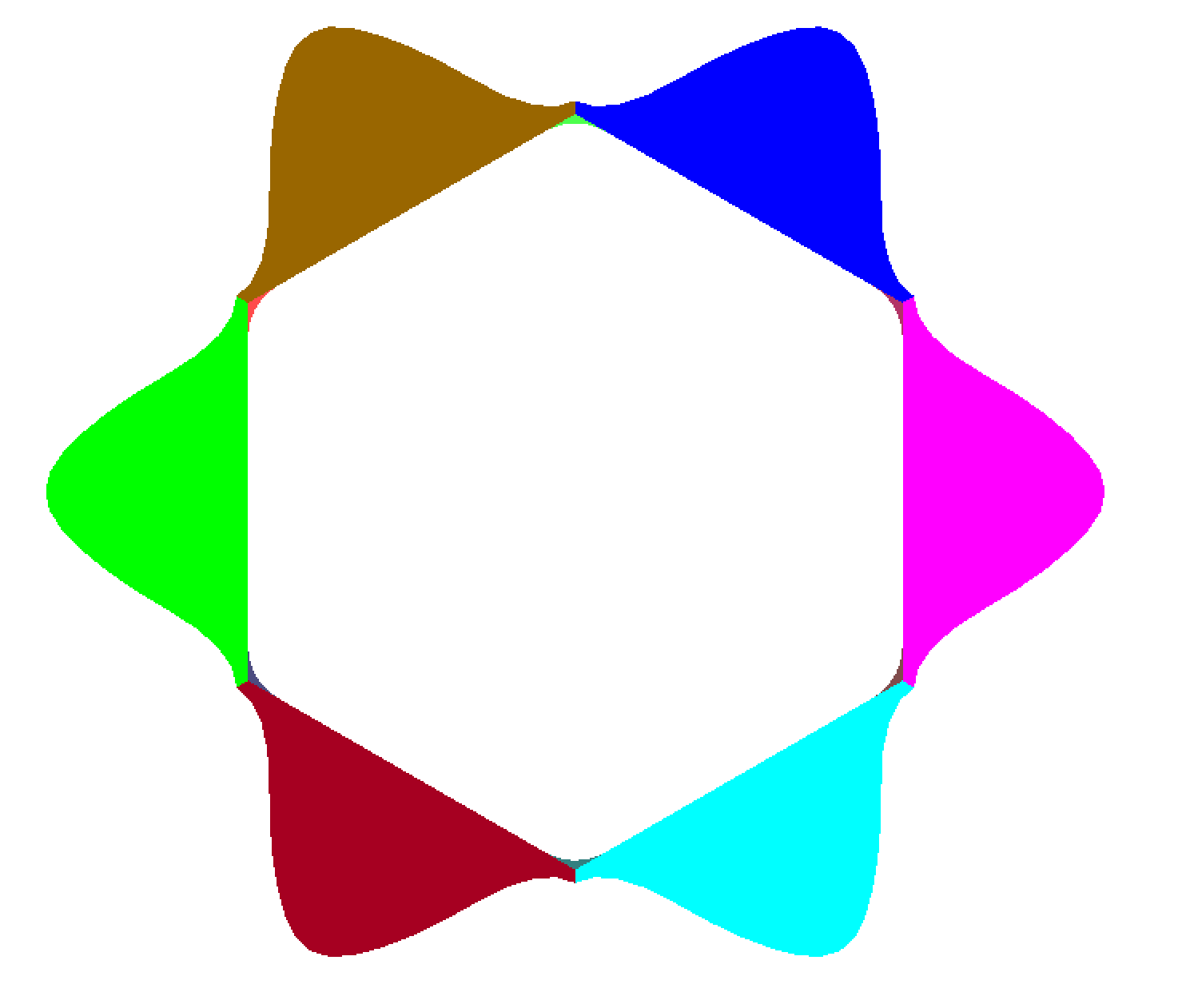}

  \end{tabular}
  \caption{Hex-can subjected to internal pressure. Deformed configuration (a)  without bending strips (b) with bending strips}.
 \label{HexcanDeformation}
\end{figure}

Here, we illustrate the method by an example representing by a Hex-can geometry subjected to internal pressure (see \Cref{StripSchematic} (Right) for the geometry of the Hex-can). The Hex-can geometry comprises of six patches joining at patch boundaries. Due to the flat geometry of each patch, kinks occur at the patch interfaces. With this configuration, the bending moment is not transmitted if C0 continuity between patches is enforced. When bending strips are added to couple the triples of control points at interfaces of each patch, the kinks do not arise (see \Cref{HexcanDeformation} (Left) for the displacement without bending strips and \Cref{HexcanDeformation} (Right) for the displacement with bending strips). As can be observed in \Cref{HexcanDeformation}, kinks do not arise, which means that bending moment is transferred between patches.

\section{FE discretization}
\label{sec:discretization}

\subsection{Discretization of the weak form}

To solve \Cref{eq:PVW} using the Newton-Raphson scheme. It shall be linearized with respect to $\Delta \mathbf{u}$:
\begin{equation}
\delta W (\bs{u}, \delta \bs{u}) + D \delta W (\bs{u}, \delta \bs{u})[\Delta \bs{u}] = \bs{0}
\end{equation}

which can be simply rewritten in a discretized form as
\begin{equation}
\delta v_a K_{ab} \Delta u_b = \delta v_a R_a
\label{eq:discretized_form}
\end{equation}

In \Cref{eq:discretized_form}, $a$ and $b$ denote the global degree of freedom (d.o.f) number of the displacement field. $R_a$ and $K_{ab}$ are the derivatives of $W$ and $R_a$, respectively:
\begin{equation}
R_a = \frac{\partial W}{\partial u_a} \qquad \qquad K_{ab} = \frac{\partial R_a}{\partial u_b} = \frac{\partial^2 W}{\partial u_a \partial u_b}
\label{eq:TangentStiffness}
\end{equation}

The residual forces vector is the difference of the internal forces vector and the external forces vector:
\begin{equation}
R_a  = R^{int}_a - R^{ext}_a,
\end{equation}

in which
\begin{align}
R^{int}_{a} &= \int\limits_{A} \left( n_{\alpha \beta}  \frac {\partial \varepsilon_{\alpha \beta}}{\partial u_a} + m_{\alpha \beta}  \frac {\partial \kappa_{\alpha \beta}}{\partial u_a} \right) \, dA \label{InternalForces} \\
R^{ext}_{a} &=  \int\limits_{A} \frac{\partial u_i}{\partial u_a} t_i \, dA = \int\limits_{A} N_I t_i \, dA
\end{align}

The tensors $n_{\alpha \beta}$ and $m_{\alpha \beta}$ in \Cref{InternalForces} denote the stress resultants, which are computed by integrating the second Piola-Kircchoff stress tensor $S^{\alpha \beta}$ over the thickness direction $\xi^3$
\begin{equation}
n_{\alpha \beta} = \int\limits_{-t/2}^{t/2} S^{\alpha \beta} \, d \xi^3 \qquad m_{\alpha \beta} = \int\limits_{-t/2}^{t/2} S^{\alpha \beta} \, \xi^3 \, d \xi^3
\end{equation}

The tangential stiffness matrix in \Cref{eq:TangentStiffness} can be split into material and geometric parts as

\begin{equation}
K_{ab} = K_{ab}^{geo} + K_{ab}^{mat}
\end{equation}

in which
\begin{align}
K_{ab}^{mat} &= \int\limits_{A} \left( \frac{\partial n_{\alpha \beta}}{\partial u_b}\frac{\partial \varepsilon_{\alpha \beta} }{\partial u_a } + \frac{\partial m_{\alpha \beta}}{\partial u_b} \frac{\partial \kappa_{\alpha \beta}}{\partial u_a  } \right) dA   \\
K_{ab}^{geo} &= \int\limits_{A} \left( n_{\alpha \beta} \frac{\partial^2 \varepsilon_{\alpha \beta} }{\partial u_a \partial u_b} + m_{\alpha \beta}\frac{\partial^2 \kappa_{\alpha \beta} }{\partial u_a \partial u_b} \right) dA \label{eq:k_geo}
\end{align}

It is noted that, the case of displacement-dependent load $t_i = t_i(u_i)$, which leads to the dependency of tangential stiffness to the internal forces, is not considered. The derivatives of the stress resultants w.r.t nodal displacements take the form of
\begin{align}
\frac{\partial n_{\alpha \beta}}{\partial u_b} &= \left( \int\limits_{-t/2}^{t/2} \hat{\mathbbm{C}}^{\alpha \beta\gamma \delta} d \xi^3 \right) \frac{ \partial \varepsilon_{\gamma \delta}}{\partial u_b} +  \left( \int\limits_{-t/2}^{t/2} \hat{\mathbbm{C}}^{\alpha \beta\gamma \delta} \xi^3 d \xi^3 \right) \frac{ \partial \kappa_{\gamma \delta}}{\partial u_b} \\
\frac{\partial m_{\alpha \beta}}{\partial u_b} &= \left( \int\limits_{-t/2}^{t/2} \hat{\mathbbm{C}}^{\alpha \beta\gamma \delta} \xi^3 d \xi^3 \right) \frac{ \partial \varepsilon_{\gamma \delta}}{\partial u_b} +  \left( \int\limits_{-t/2}^{t/2} \hat{\mathbbm{C}}^{\alpha \beta\gamma \delta} (\xi^3)^2 d \xi^3 \right) \frac{ \partial \kappa_{\gamma \delta}}{\partial u_b}
\end{align}

For more details on the derivatives of strain tensor with respect to nodal displacement, the reader is referred to \citep{Kiendl2015,Kiendl2011}. Note that the four-order tensor $\hat{\mathbbm{C}}^{\alpha \beta\gamma \delta}$ in the equations above is expressed in \Cref{app:tangent_moduli} which is then modified following Eq.(\ref{CondensedTensor}).

\subsection{B\'ezier decomposition}

To evaluate basis function \Cref{NURBSBasisFunctions}, the full knot vector is required. In addition, to interpolate the displacement within an element, the displacement at other control points in the patch is necessary. To alleviate this issue and maintain the local characteristic of the finite element, the B\'ezier decomposition concept is proposed, which bridges the gap between isogeometric method and standard finite element method \citep{Borden2011}. Taking into account that not all B-splines basis function is nonzero on a knot span, the basis functions $\mathbf{N}$, which is not vanished on a knot span, can be represented by a linear combination of B\'ezier basis functions $\mathbf{B}$ constructed on knot vector $\{ \underbrace{0,\hdots,0}_{p+1 \, \text{times}}, \underbrace{1,\hdots,1}_{p+1 \, \text{times}} \}$, in short:
\begin{equation}
\mathbf{N} = \mathcal{C} \mathbf{B}
\end{equation}

\begin{figure}[h!]
\centering
\includegraphics[scale=0.1]{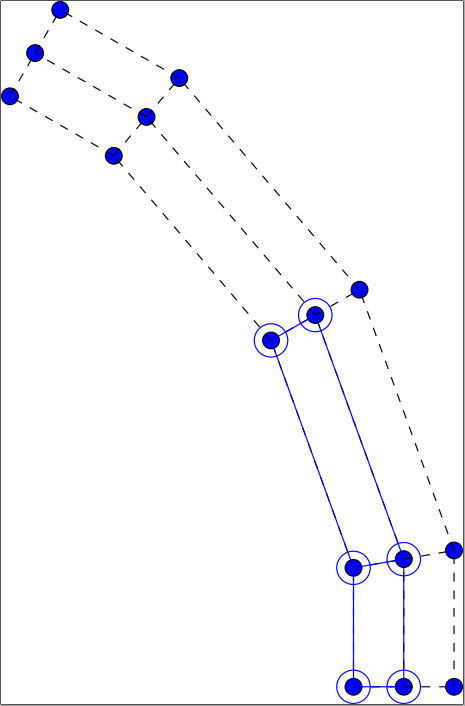}
\includegraphics[scale=0.1]{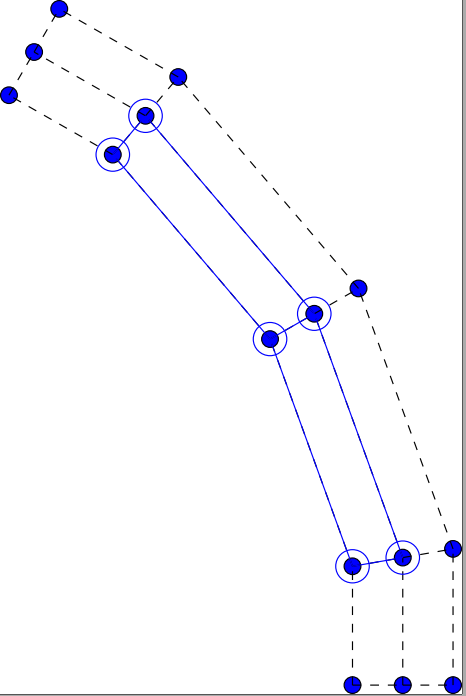}
\includegraphics[scale=0.1]{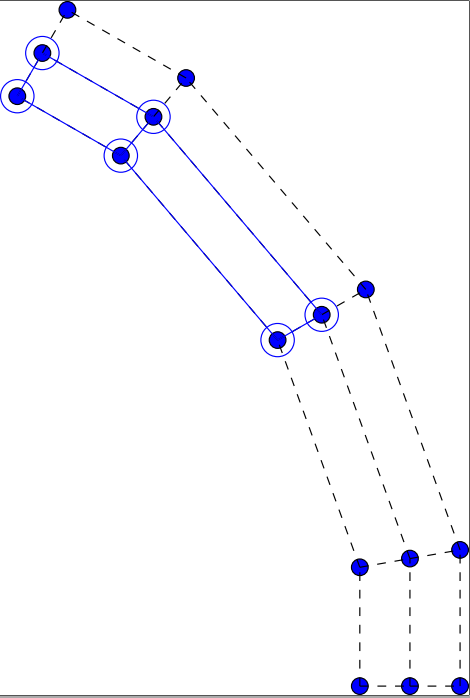}
\includegraphics[scale=0.1]{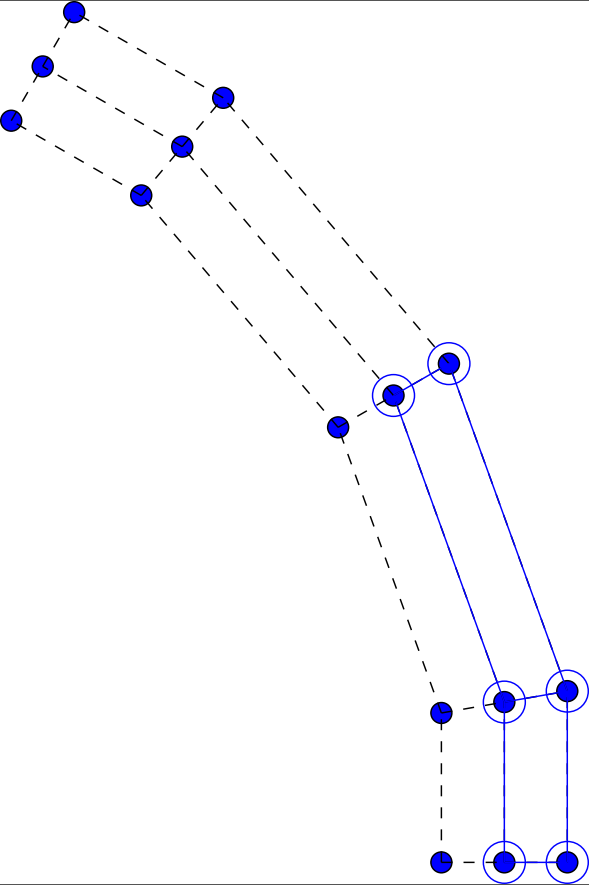}
\includegraphics[scale=0.1]{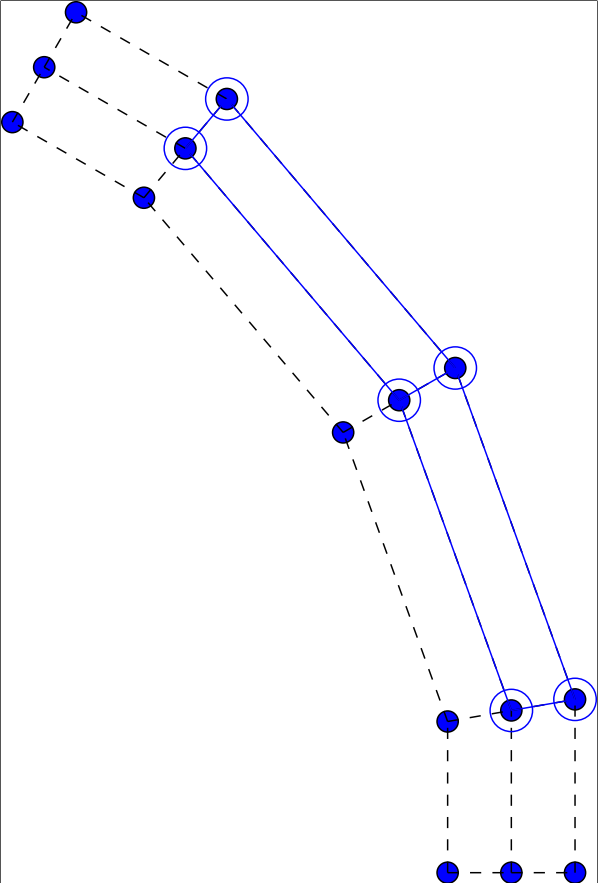}
\includegraphics[scale=0.1]{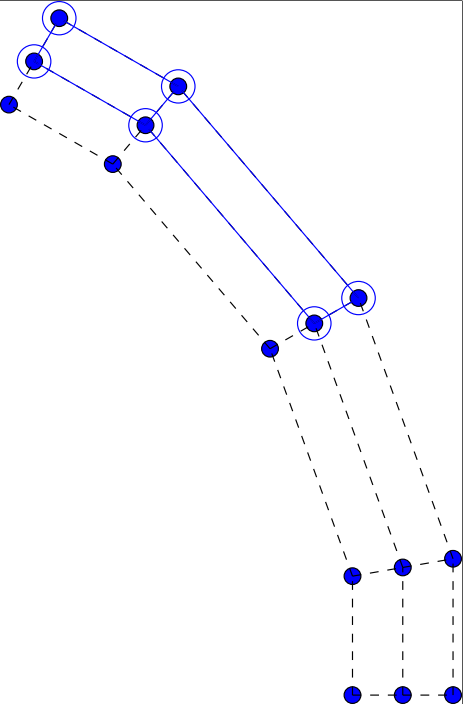}
\caption{B\'ezier elements over the isogeometric control mesh; \protect\tikz{\protect\draw[fill=blue] (0,0) circle (.7mm);} : control point; \protect\tikz{\protect\draw[fill=blue] (0,0) circle (.7mm);\protect\draw[] (0,0) circle (1.4mm);} : control point supported by B\'ezier element.}
\label{fig:bezier_elements}
\end{figure}

The coefficient matrix $\mathcal{C}$ is called B\'ezier extraction operator and depends only on the knot vectors, hence it is considered as constant. For more details on how to compute $\mathcal{C}$, the reader is referred to \citep{Borden2011}. 

Based on the B\'ezier decomposition concept, the B\'ezier finite element can be constructed, which contains only the control points associated with non-vanished basis function on the knot span. An illustration of the B\'ezier finite element for 2D case can be found in \Cref{fig:bezier_elements}.

To evaluate the second derivatives terms, i.e. \Cref{eq:k_geo}, the second derivatives of the basis functions are required. \Cref{app:bezier_second_derivatives} presents a method to take the second derivatives of NURBS basis functions evaluated through B\'ezier extraction operator and B\'ezier basis functions.

\section {Numerical examples}
\label{sec:examples}

\begin{figure}[htb]
		\color{blue}
	\centering
	\includegraphics[scale=0.5]{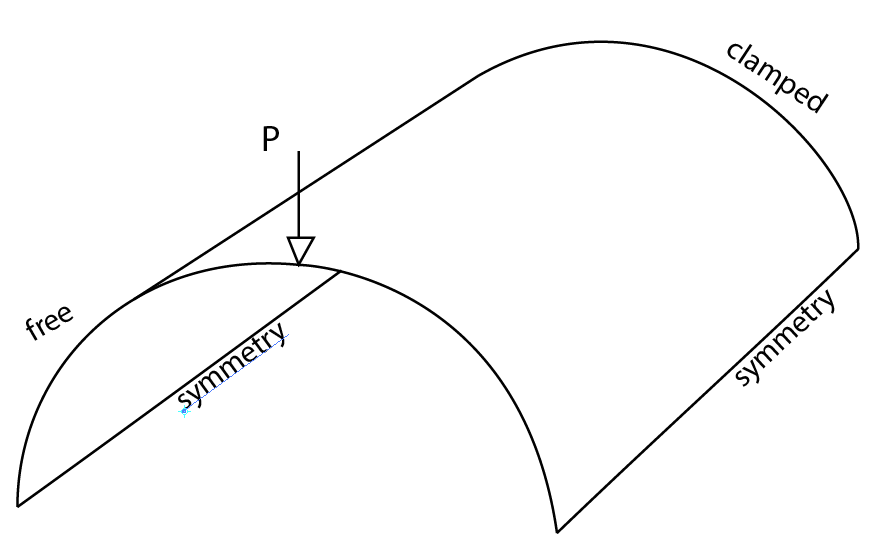}
	\caption{\color{blue}{Numerical example 1: Geometry and boundary condition. }}
	\label{PinchSemiCylinder}
\end{figure}

In this section, four selected examples are performed to analyze hyperelastic and elastoplastic behaviours of multi-patch thin shell structures. Comparisons with the reference results are made. Arc-length control method following the approach in \citep{Crisfield1981} is adopted to capture the structural response in the last two examples.
\textcolor{blue}{
\subsection{Numerical example 1: Pinched semi-cylindrical shell}
For the first example, a semi-cylidrical shell is pinched by an end force at the middle of one end and fully clamped in the other direction. The point load is applied incrementally until the maximum load $P=2000.0$ is reached. Figure. \ref{PinchSemiCylinder} shows the problem setup and boundary conditions. The length of the cylinder is $L=3.048$ and the radius is $R=1.016$. St. Vennant-Kirchhoff material model is used with $E=2.0685 \times 10^7$ and $\rho = 0.3$.}
	
\begin{figure}[htb]
		
	\centering
	\includegraphics[scale=0.5]{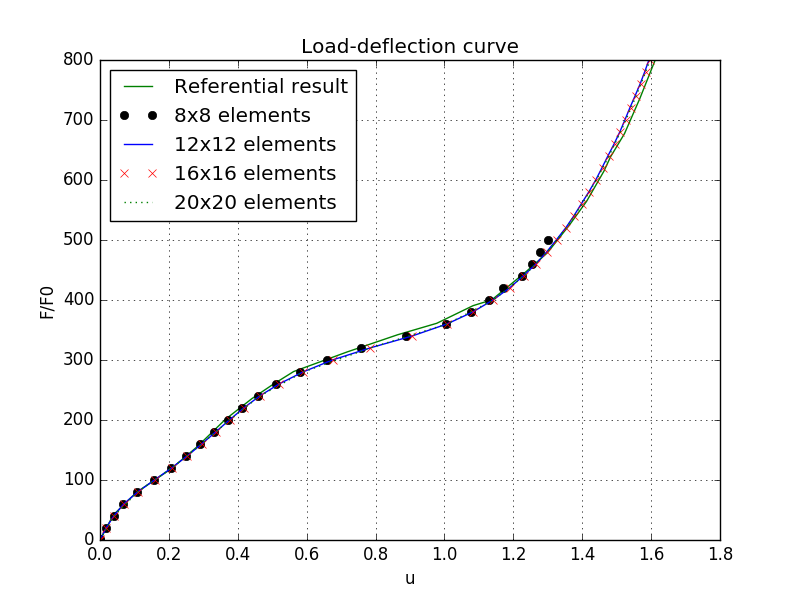}
	\caption{\color{blue}{Numerical example 1: Investigation of mesh sensitivity with respect to the load-deflection curve.} }
	\label{PinchSemiCylinderLD}
\end{figure}

\begin{figure}[htb]
	\color{blue}
	\centering
	\includegraphics[scale=0.5]{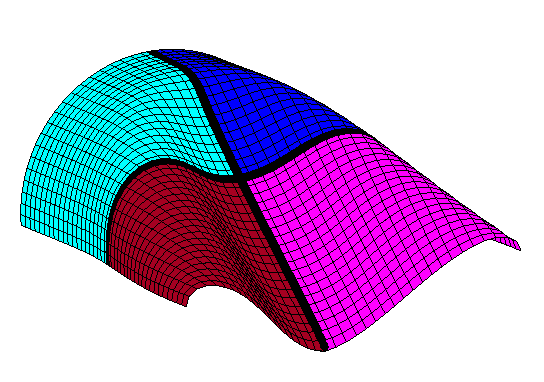}
	\caption{Numerical example 1: Deformed geometry at the maximum load. }
	\label{DeformedClampedCyl}
\end{figure}

\textcolor{blue}{Figure.\ref{PinchSemiCylinderLD} reveals obtained results on four different mesh sizes with cubic NURBS elements, indicating that the computational procedure fails at $P=500.0$ with $8 \times 8$ elements and accurate results are obtained with the mesh sizes of $12 \times 12$, $16 \times 16$ and $20 \times 20$. The deformed geometry of the cylinder under the maximum load is revealed in Fig.\ref{DeformedClampedCyl}. Table.\ref{IterationCompar} reports the total number of iterations required for the computation, in which the number for our approach that requires 3 to 9 iterations at each load step is less than the one for other shell elements \cite{Sze2004} \cite{Brank1995} using finite element methods . }

\begin{table}[h]
	\centering
	\color{blue}
	\begin{tabular} {l c c c c}
		\hline
					 		  & 16x16  & 32x32 S4R elements & 16x16  &16x16 elements  \\
					 		  & current elements &  \cite{Sze2004}   & current elements & ,B. Brank et al. \cite{Brank1995}\\
					 		  \hline
		Incremental load      & 40  			&  33		&80		& 80  \\
		Number of iterations  & 124				& 184		&76		& 110 \\
		\hline

	\end{tabular}
	\caption{\color{blue}{Numerical example 1: Lists of the numbers of iterations and load increments.} }
\label{IterationCompar}
\end{table}
\subsection{Numerical example 2: Pinching of cylinder}

\begin{figure}[htb]
\centering
  \includegraphics[scale=0.7]{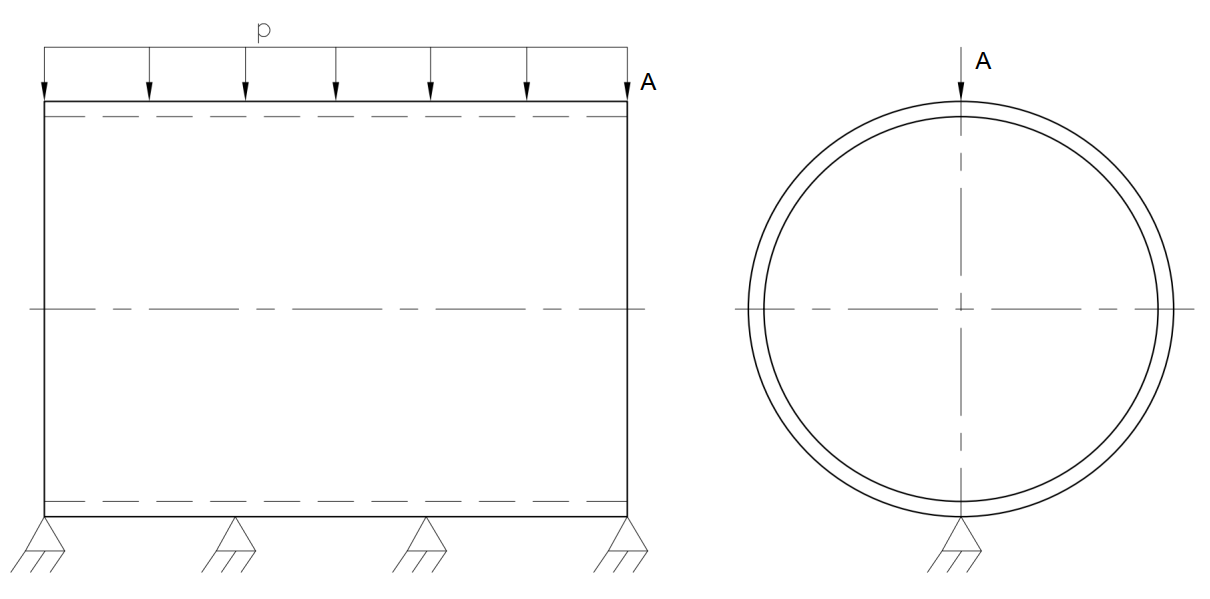}
\caption{Numerical example 2: Geometry and boundary condition }
\label{PinchCylinderGeometry}
\end{figure}

In the second example, the proposed computational model is validated against finite strain hyperelasticity model. Following that, a tube with geometry and boundary conditions as illustrated in \Cref{CylinderGeometry} is analyzed. The geometric parameters are chosen as length $L=30 cm$, radius $R=9 cm$, and thickness $t=0.2 cm$. The bottom of the cylinder is fixed and a uniformly distributed line load ($p$) is applied on top (see \Cref{PinchCylinderGeometry}). The constitutive law is Neo-Hookean with the strain energy is defined in \Cref{eq:neo_hookean}. The material parameters are $\mu = 60 kN/mm^2$ and $\lambda = 240 kN/mm^2$. The geometry of the cylinder consists of 4 patches, and is discretized using 552 cubic B\'ezier elements.

\begin{equation}
W = \dfrac{\mu}{2} \left( \text{tr} \bs{C} - 3 \ \right) - \mu \ln{\sqrt{\det{\bs{C}}}} + \dfrac{\lambda}{4} \left( \det{\bs{C}} - 1 - 2 \ln{\sqrt{\det{\bs{C}}}} \right)
\label{eq:neo_hookean}
\end{equation}

\begin{figure}[htb]
\centering
\includegraphics[scale=0.35]{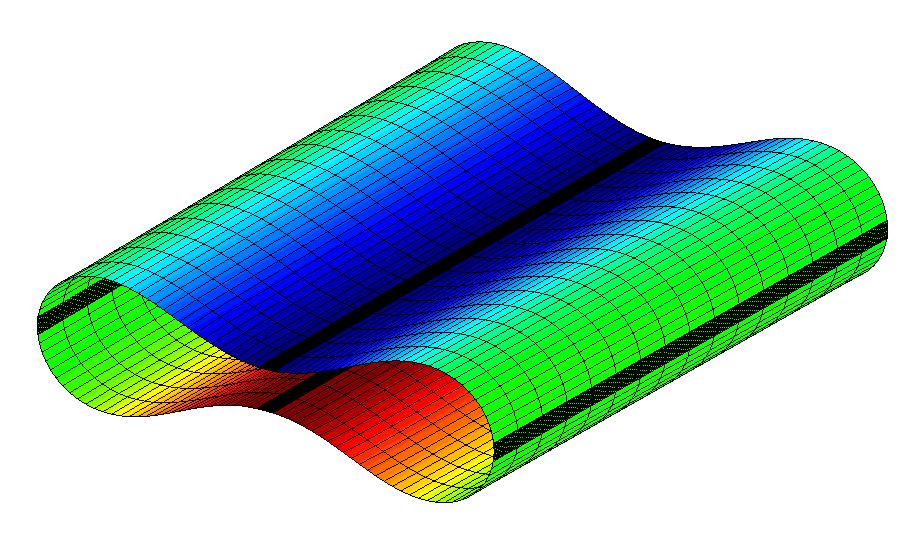}
\includegraphics[scale=0.3]{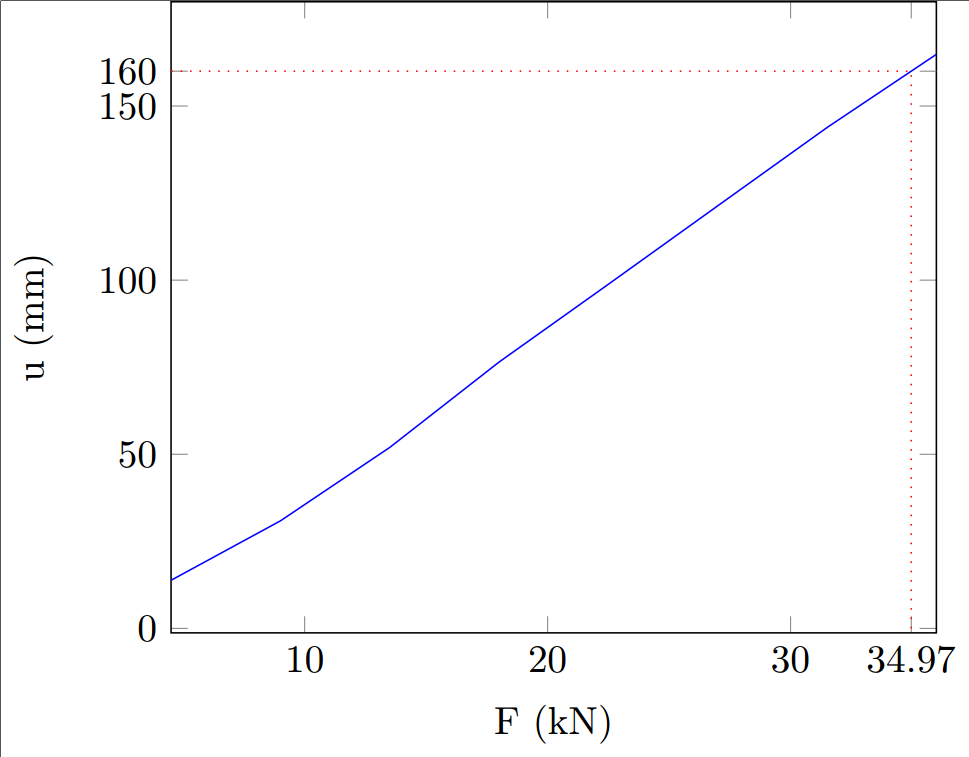}
\caption{Numerical example 2: Left) Contour plot of the last deformed configuration; Right) Load-displacement curve measured on top.}
\label{DeformedPinchingCylinder}
\end{figure}

The load $p$ is increased equally in 8 steps. \Cref{DispList} lists computed displacements at point A corresponding to applied load at each  step. From \Cref{DeformedPinchingCylinder} (Right), the structural load when the vertical deflection is $u=160 mm$ can be approximated as $F\approx 34.965 kN$, which is in good agreement with the results from literature ranging between $34.59 kN$ and $35.47 kN$ \citep{Schwarze2011}. A contour plot of the last deformed configuration is depicted in \Cref{DeformedPinchingCylinder} (Left).

\begin{table}[h]
\centering
  \begin{tabular}{c c c c c c c c c}
  \hline
 Steps        &1 & 2 & 3 &4 & 5 & 6 & 7 & 8 \\
  \hline
  $F (kN)    $ & $4.5$      & $9.0       $ & $13.5   $ & $18.0   $ & $22.5     $ & $27.0     $  & $31.5  $      & $36.0 $ \\
  $u (mm) $ & $13.817$ & $30.872  $ & $51.954$ & $76.473$ & $98.818  $ & $121.337$  & $143.815  $ & $164.83$ \\
  \hline
  \end{tabular}
\caption{Numerical example 2: Resulting displacements corresponding to applied load at each step.}
\label{DispList}
\end{table}

\subsection{Numerical example 3: Scordelis-Lo roof}
\begin{figure}[htb]
\centering
  \includegraphics[scale=0.25]{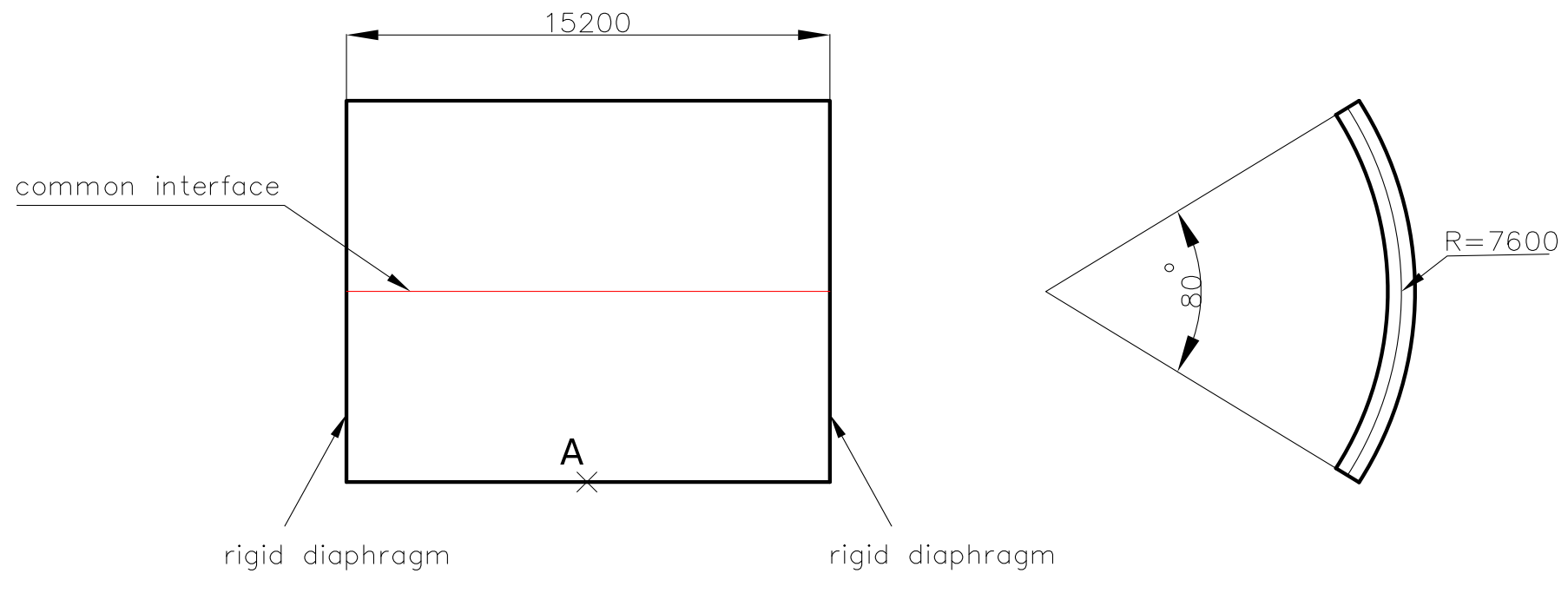}
\caption{Numerical example 3: Geometry and boundary conditions. It is noted that, the rigid diaphragms are used to support the two ends of the roof.}
\label{ScordelisGeometry}
\end{figure}

In the next example, the failure analysis of the Scordelis-Lo roof, accounting for fully non-linear thin shell analysis, is considered. \Cref{ScordelisGeometry} illustrates  boundary conditions and the setup of the geometry of the problem which contains 2 patches connected at the middle line on the roof top.The material parameters are chosen as elastic modulus $E = 2.1 \times 10^4 N/mm^2$,  Poisson's ratio $\rho=0$ and the hardening function $\kappa(\alpha) = 4.2 N/mm^2$. Because the hardening modulus is vanished, the perfect plasticity condition is obtained. The structure is applied with gravity load. A mesh comprises of $3969$ quadratic B\'ezier elements is used to discretize the whole domain. Arc-length control method is employed to capture the structural response with the reference value $f_0 = 4 \times 10^{-3} N/mm^{2}$.

\begin{figure}[htb!]
\centering
  \begin{tabular}{@{}cccc@{}}
    \includegraphics[width=0.35\textwidth]{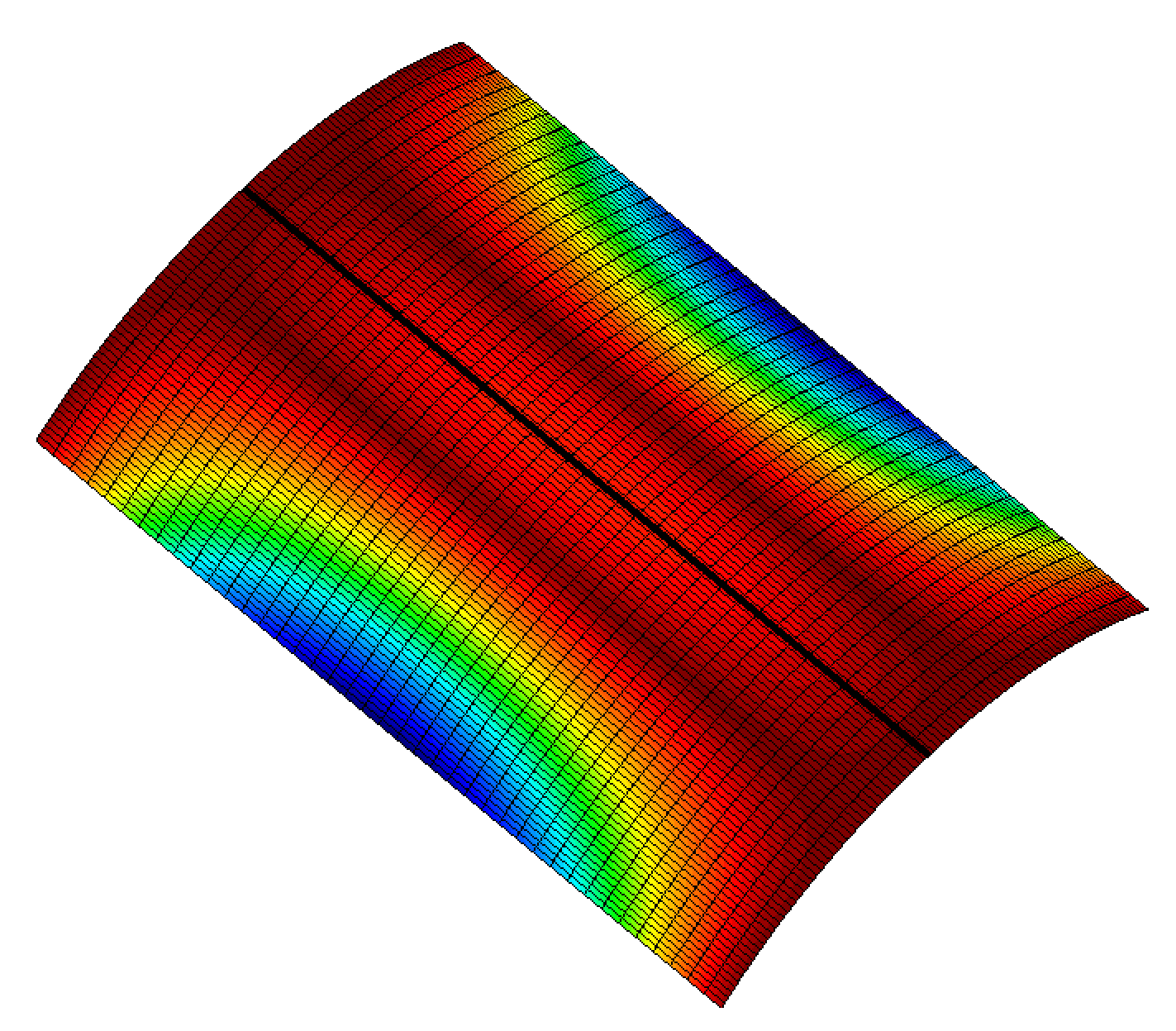}
    \includegraphics[width=0.35\textwidth]{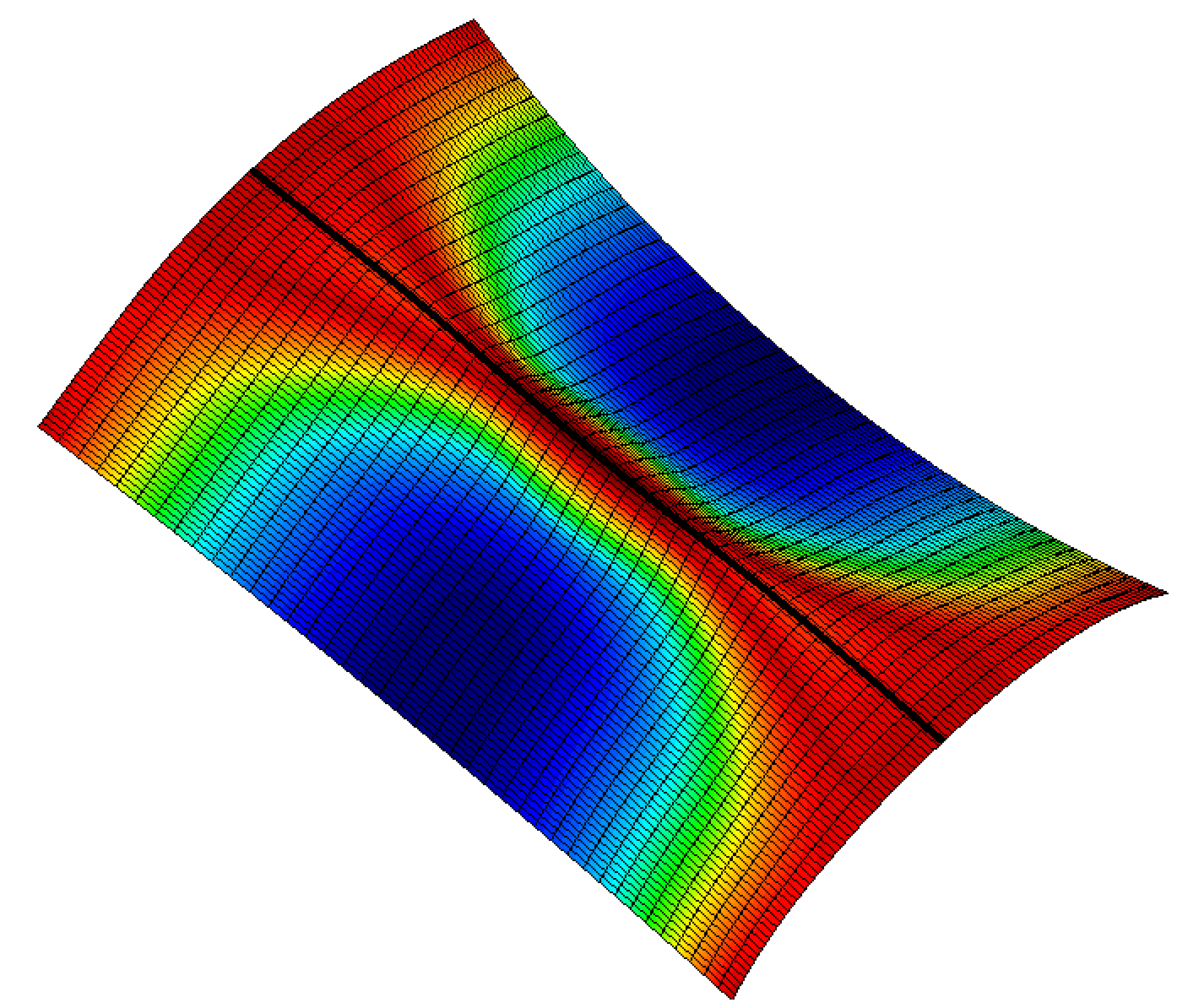}&  \\[0.3cm]
    \includegraphics[width=0.35\textwidth]{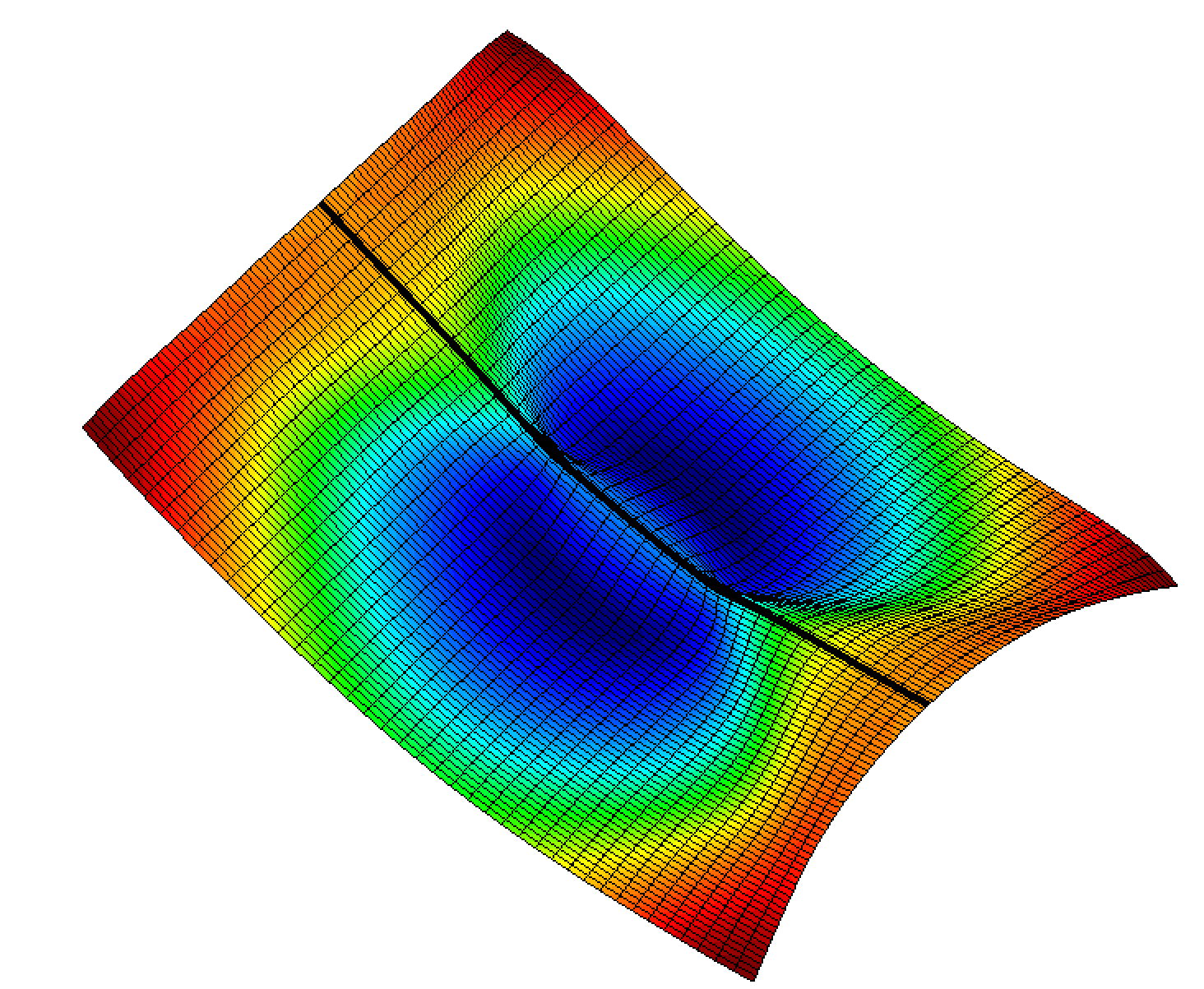}
    \includegraphics[width=0.35\textwidth]{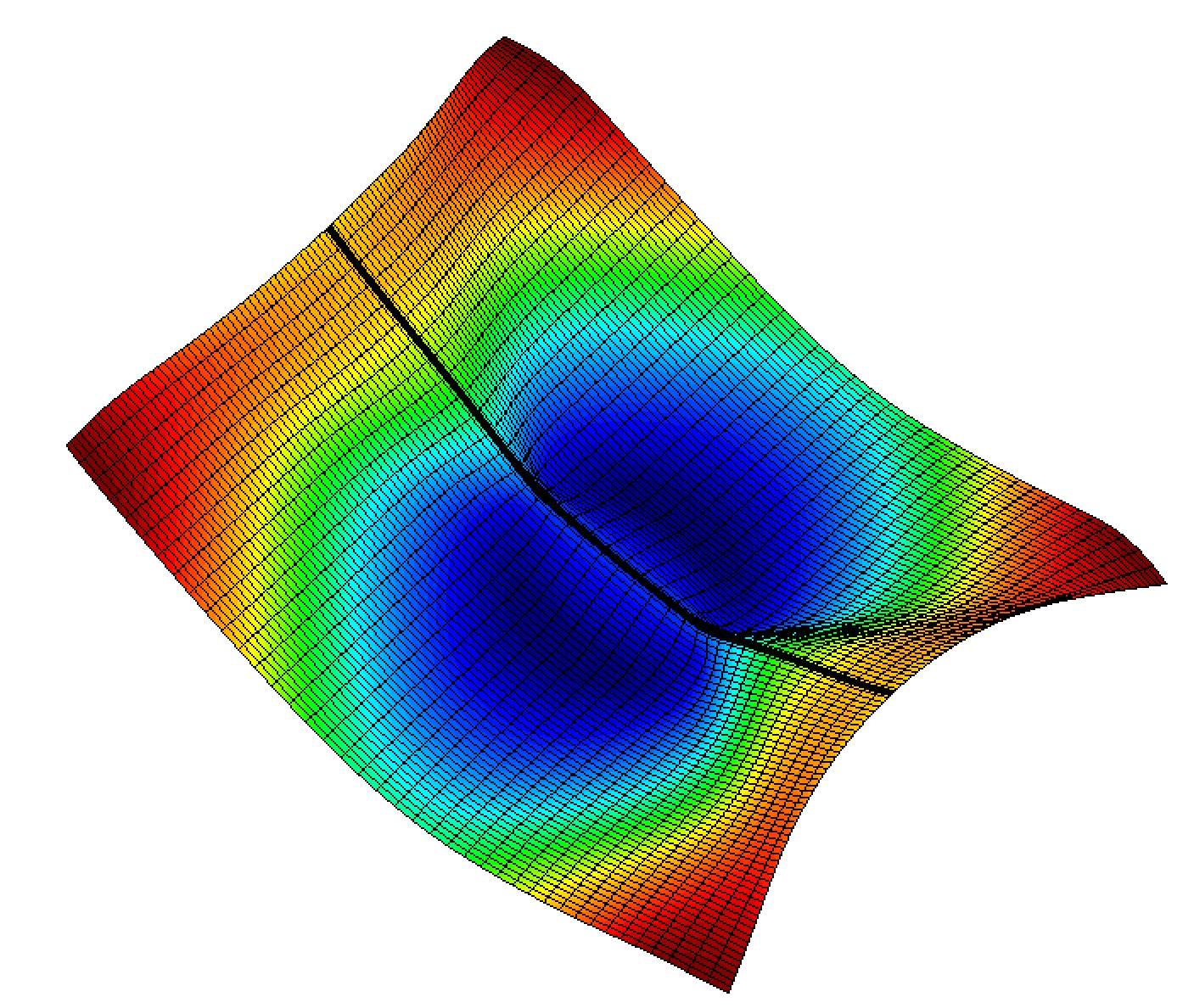}    &  \\[0.3cm]
  \end{tabular}
  \caption{Numerical example 3: Deformation at different stages.}
 \label{ScodelisDeformation}
\end{figure}

\Cref{ScodelisDeformation} shows deformed configurations at different loading stages. As can be observed from \Cref{ScodelisDeformation}, the localized failure mode is characterized by the appearance of the plastic hinge along the axial line on the roof top.

\begin{figure}[htb!]
\centering
\begin{tikzpicture}
  \begin{axis}[
  ymin = 0,
  ymax = 1.5,
  xmin = 0,
  xmax = 2500,
  xlabel = u (mm),
  ylabel = $f/f_0$,
  xtick={0, 500, 1000, 1500, 2000, 2500},
  legend pos=south west,
  legend style={at={(0.6,0.65)}},
  ]
  \addplot[mark=none,red,dotted,line width=0.3mm] table[x expr=1e3*\thisrowno{0},y expr=\thisrowno{1}] {graph/scordelis_lo_roof_brank.txt};
  \addlegendentry{Reference \citep{Brank1997}}
  \addplot[mark=none,cyan,dashdotted,line width=0.3mm] table[x expr=\thisrowno{0},y expr=\thisrowno{1}] {graph/scordelis_lo_roof_ambati.txt};
  \addlegendentry{Reference \citep{Ambati2018}}
  \addplot[mark=none,blue,smooth,line width=0.3mm] table[x expr=\thisrowno{0},y expr=\thisrowno{1}] {graph/scordelis_lo_roof_dg.txt};
  \addlegendentry{Obtained results}
  \end{axis}
\end{tikzpicture}
\caption{Numerical example 3: Load-deflection curve at point A.}
\label{ScodelisLoadDeflection}
\end{figure}
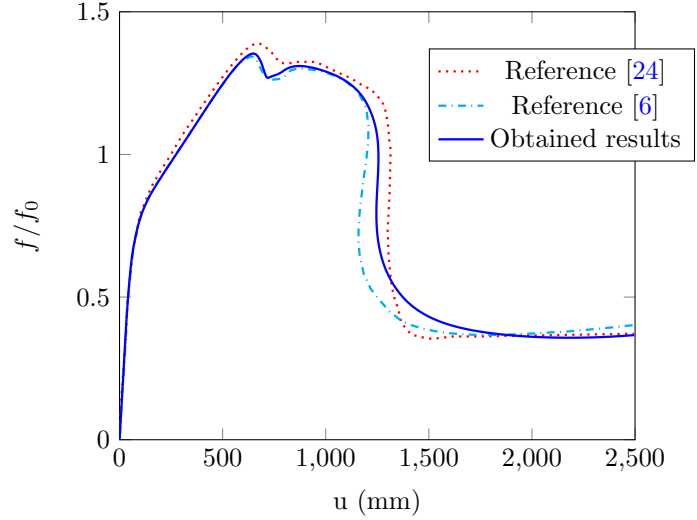

The corresponding load-deflection curve at center point of the side (denoted A) is plotted in \Cref{ScodelisLoadDeflection}. It can be observed that the  obtained results agree well with the reference\citep{Ambati2018,Brank1997}.

\subsection{Numerical example 4: Pinching of cylinder with large strains}
\begin{figure}[h]
\centering
  \includegraphics[scale=0.25]{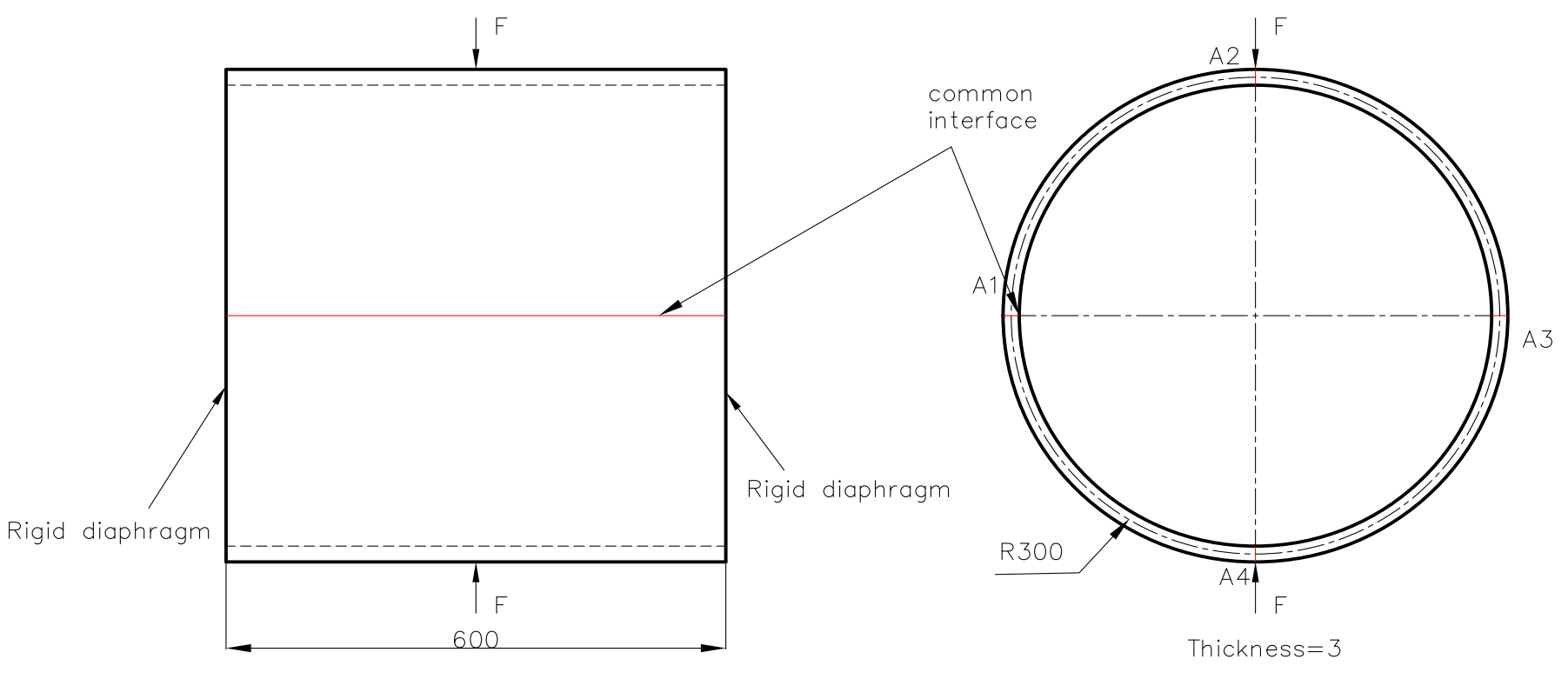}
\caption{Numerical example 4: Geometry and boundary conditions. (A1, A2, A3 and A4 are marked in red.)}
\label{CylinderGeometry}
\end{figure}

\begin{figure}[!htb]
\centering
  \begin{tabular}{@{}cccc@{}}
    \includegraphics[width=0.35\textwidth]{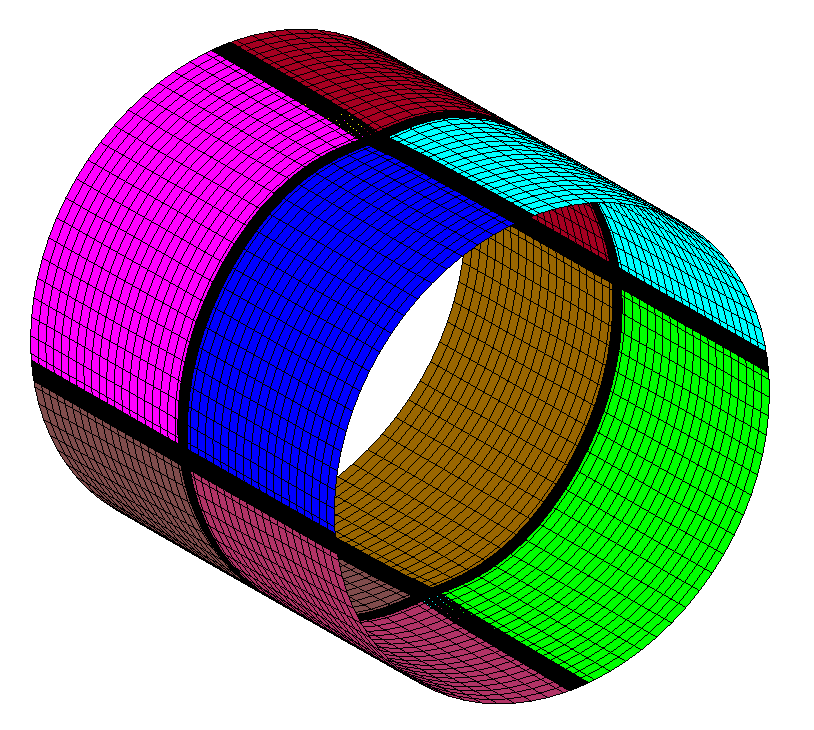}
    \includegraphics[width=0.35\textwidth]{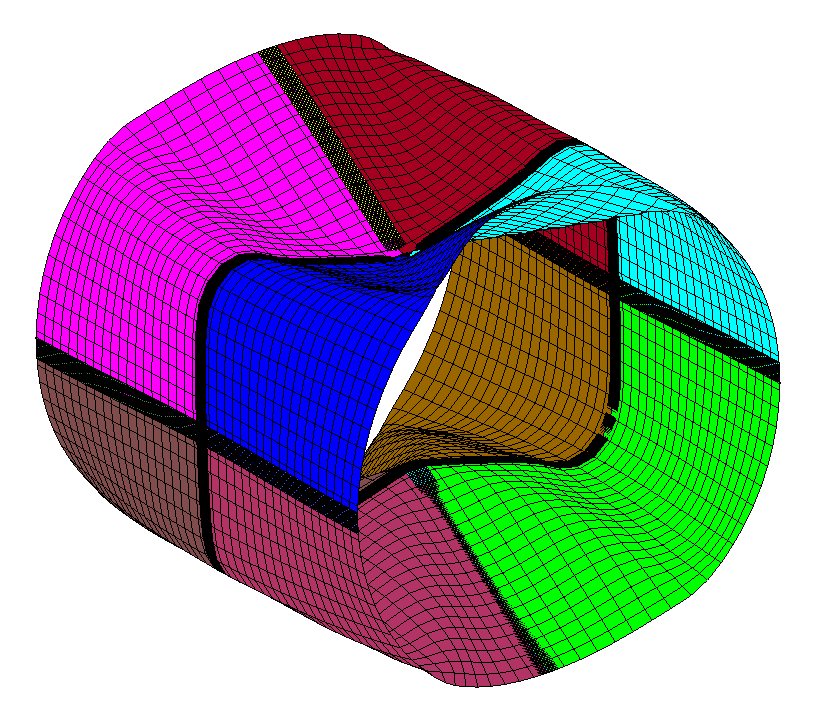} &  \\[0.2cm]
    \includegraphics[width=0.35\textwidth]{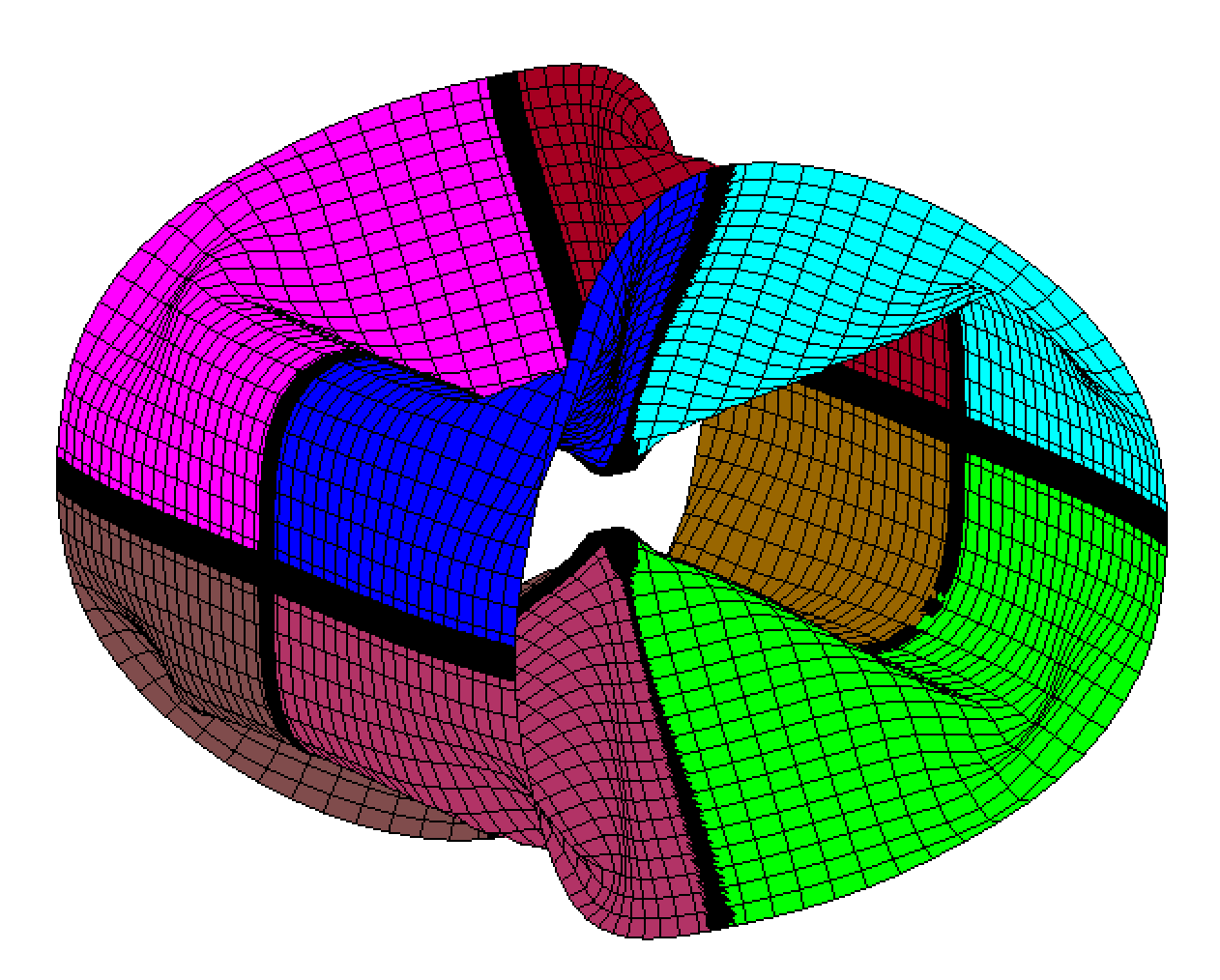}
    \includegraphics[width=0.35\textwidth]{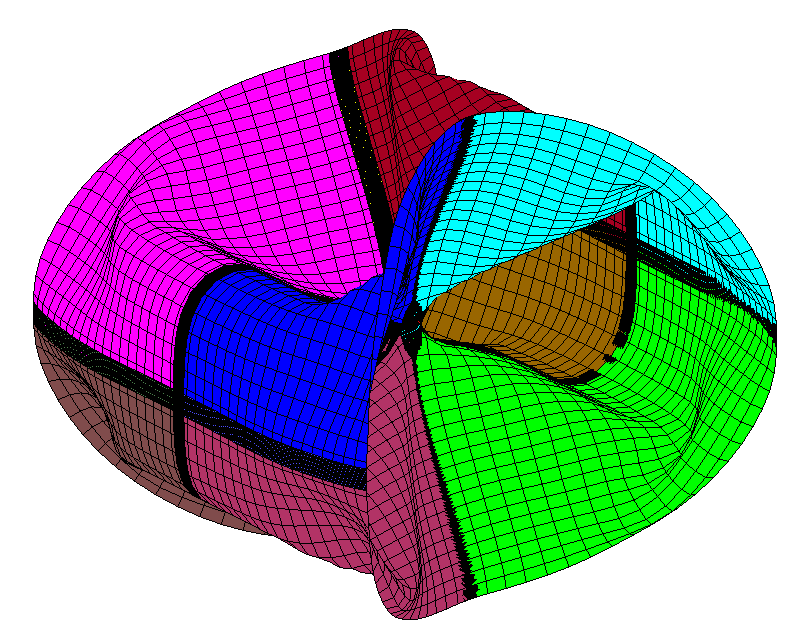}  &  \\[0.2cm]
  \end{tabular}
  \caption{Numerical example 4: Deformation at different stages.}
 \label{CylinderDeformation}
\end{figure}

The last example considers a cylinder supported by two rigid diaphragms and pinched by two concentrated forces at the middle of the opposite sides. \Cref{CylinderGeometry} presents the setup of boundary conditions and the geometry of the cylinder which contains four patches connected at common interfaces A1, A2, A3 and A4. The material parameters are chosen as $E =3000 $, $\rho=0.3$ and the hardening function $\kappa(\alpha) = 24.3 + 300 \alpha$. The geometry contains 8 patches as depicted in \Cref{CylinderDeformation} (Top Left). The bending strips are used on each patch interfaces to enforce the $C^1$-continuity condition. A mesh contains 3452 quadratic B\'ezier elements is used to discretize the whole domain.

\begin{figure}[!htb]
\centering
\begin{tikzpicture}
  \begin{axis}[
  ymin = 0,
  ymax = 35,
  xmin = 0,
  xmax = 300,
  xlabel = u (mm),
  ylabel = F (kN),
  xtick={0, 500, 1000, 1500, 2000, 2500},
  legend pos=north west,
  ]
%
  \addplot[mark=none,red,dotted,line width=0.3mm] table[x expr=\thisrowno{0},y expr=\thisrowno{1}*4/1000] {graph/pinched_cylinder_plasticity_areias.txt};
  \addlegendentry{Reference \citep{Pedro2016}}
  \addplot[mark=none,cyan,dashdotted,line width=0.3mm] table[x expr=\thisrowno{0},y expr=\thisrowno{1}/1000] {graph/pinched_cylinder_plasticity_ambati.txt};
  \addlegendentry{Reference \citep{Ambati2018}}
  \addplot[mark=none,blue,smooth,line width=0.3mm] table[x expr=\thisrowno{0},y expr=\thisrowno{1}/1000] {graph/pinched_cylinder_plasticity_dg.txt};
  \addlegendentry{Obtained results}
  \end{axis}
\end{tikzpicture}
\caption{Numerical example 4: Load-deflection curve.}
\label{CylinderLoadDeflection}
\end{figure}
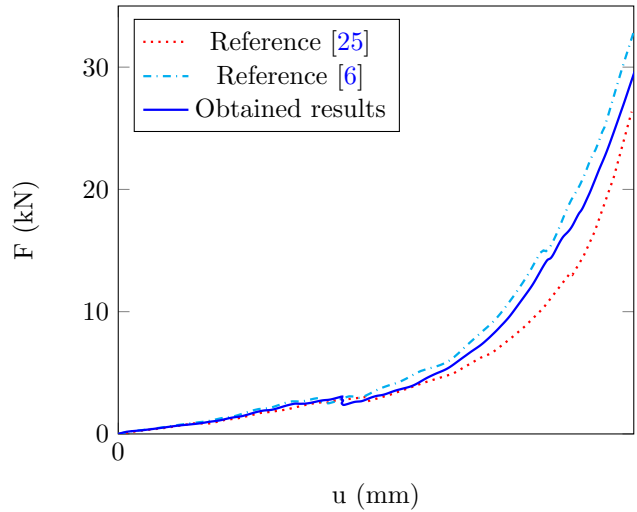

Deformed configurations at different stages of loading are illustrated in \Cref{CylinderDeformation}. The load-deflection curve representing the structural response, obtained at the loading point, is shown in \Cref{CylinderLoadDeflection}. The obtained results are in good agreement with the reference \citep{Ambati2018,Pedro2016}.


\section{Conclusion}
\label{sec:conclusion}
An efficient computational model for thin shell using isogeometric analysis is proposed, with the advantage of using bending strip method to maintain the $C^1$ continuity, thus allowing for efficient modelling using multiple NURBS patch with $C^0$ continuity at the patch boundaries. The B\'ezier decomposition method is used to retain the local characteristics of the finite element. This computational model is validated with respect to nonlinear constitutive model, including hyperelasticity and elastoplasticity. The benchmark results agree well with the references, showing that the proposed computational model is promising for efficient shell analysis.\textcolor{blue}{ The presented shell formulation is very general and flexible and can be extended to dynamics and ductile fracture motivated by published papers that are classified into two categories: discrete crack modeling \cite{Areias2005b} \cite{Song2009} \cite{Song2006} and smeared crack modeling \cite{Areias2016} \cite{Reinoso2017} .} On the other hand, in order to further advance the efficiency, the parallelization for distributed computing shall be done. The parallel algorithm shall take into account the distribution of bending strips over computing processes. That will be addressed in future research.

\section*{Acknowledgements}

The first author would like to acknowledge the financial support via RISE project BESTOFRAC 734370 for this work. The research performed by Hoang-Giang Bui and Günther Meschke were conducted in the framework of the Collaborative Research Project SFB 837 "Interaction Modeling in Mechanized Tunneling", financed by the German Research Foundation (DFG). The authors would like to thank the DFG for the support of this project.\\

\begin{appendices}

\crefalias{section}{appsec}
\crefalias{subsection}{appsec}

\section{Second derivatives of NURBS basis functions based on B\'ezier decomposition}
\label{app:bezier_second_derivatives}

The NURBS basis function, which is non-vanished on the knot span, is defined by
\begin{equation}
\mathbf{R} = \dfrac{\mathcal{W} \mathcal{C} \mathbf{B}}{\langle \mathbf{W}^b, \mathbf{B} \rangle}
\label{eq:nurbs_bf}
\end{equation}

In \Cref{eq:nurbs_bf}, the operator $\langle \cdot, \cdot \rangle$ denotes the dot product of 2 vectors. $\mathcal{W}$ is the diagonal matrix with entries as the corresponding weight of the non-vanished basis function. $\mathbf{W} = \text{diag} \{\mathcal{W}\}$ and $\mathbf{W}^b = \mathcal{C}^T \mathbf{W}$.

The first derivatives of \Cref{eq:nurbs_bf} can be computed as
\begin{equation}
\dfrac{\partial \mathbf{R}}{\partial \xi^\alpha} = \mathcal{W} \mathcal{C} \left[ \dfrac{ \nicefrac{\partial \mathbf{B}}{\partial \xi^\alpha}  }{\langle \mathbf{W}^b, \mathbf{B} \rangle} - \dfrac{ \langle \mathbf{W}^b, \nicefrac{\partial \mathbf{B}}{\partial \xi^\alpha} \rangle \mathbf{B} }{\langle \mathbf{W}^b, \mathbf{B} \rangle^2} \right]
\label{eq:nurbs_bf_1st_der}
\end{equation}

Taking the derivatives of \Cref{eq:nurbs_bf_1st_der} with respect to $\xi^\alpha$ and $\xi^\beta$ gives
\begin{equation}
\dfrac{\partial^2 \mathbf{R}}{\partial (\xi^\alpha)^2} =
\mathcal{W} \mathcal{C} \left[ \dfrac{ \nicefrac{\partial^2 \mathbf{B}}{\partial (\xi^\alpha)^2} }{\langle \mathbf{W}^b, \mathbf{B} \rangle}
- 2 \dfrac{ \langle \mathbf{W}^b, \nicefrac{\partial \mathbf{B}}{\partial \xi^\alpha} \rangle \nicefrac{\partial \mathbf{B}}{\partial \xi^\alpha} }{\langle \mathbf{W}^b, \mathbf{B} \rangle^2}
+ 2 \dfrac{ \langle \mathbf{W}^b, \nicefrac{\partial \mathbf{B}}{\partial \xi^\alpha}^2 \rangle \mathbf{B} }{\langle \mathbf{W}^b, \mathbf{B} \rangle^3}
- \dfrac{ \langle \mathbf{W}^b, \nicefrac{\partial^2 \mathbf{B}}{\partial (\xi^\alpha)^2} \rangle \mathbf{B} }{\langle \mathbf{W}^b, \mathbf{B} \rangle^2}\right]
\end{equation}

And
\begin{equation}
\begin{split}
\dfrac{\partial^2 \mathbf{R}}{\partial \xi^\alpha \partial \xi^\beta} = \mathcal{W} \mathcal{C} \bigg[ &\dfrac{ \nicefrac{\partial^2 \mathbf{B}}{\partial \xi^\alpha \partial \xi^\beta} }{\langle \mathbf{W}^b, \mathbf{B} \rangle}
- \dfrac{ \langle \mathbf{W}^b, \nicefrac{\partial \mathbf{B}}{\partial \xi^\alpha} \rangle \nicefrac{\partial \mathbf{B}}{\partial \xi^\beta} }{\langle \mathbf{W}^b, \mathbf{B} \rangle^2}
- \dfrac{ \langle \mathbf{W}^b, \nicefrac{\partial \mathbf{B}}{\partial \xi^\beta} \rangle \nicefrac{\partial \mathbf{B}}{\partial \xi^\alpha} }{\langle \mathbf{W}^b, \mathbf{B} \rangle^2} \\
& + 2 \dfrac{ \langle \mathbf{W}^b, \nicefrac{\partial \mathbf{B}}{\partial \xi^\alpha} \rangle \langle \mathbf{W}^b, \nicefrac{\partial \mathbf{B}}{\partial \xi^\beta} \rangle \mathbf{B} }{\langle \mathbf{W}^b, \mathbf{B} \rangle^3}
- \dfrac{ \langle \mathbf{W}^b, \nicefrac{\partial^2 \mathbf{B}}{\partial \xi^\alpha \partial \xi^\beta} \rangle \mathbf{B} }{\langle \mathbf{W}^b, \mathbf{B} \rangle^2} \bigg]
\end{split}
\end{equation}



\section{Consistent elasto-plastic tangent moduli}
\label{app:tangent_moduli}

\begin{table}
  \begin{center}
    \begin{tabular}{l} 
       \toprule[2 pt]
       Scaling factors:\\
\hspace*{5mm}$ \bar{G} = \frac{1}{3} G \text{tr}(\bar{\bs{b}}^e) $  \\
\hspace*{5mm}$ \Xi_0 = 1+ \frac{R'}{3 \bar{G}}$  \\
\hspace*{5mm}$ \Xi_1 = [1 - \frac{1}{\Xi_0}] \frac{2}{3} \frac{||\bs{\tau}_{dev,n+1}^{trial}||}{\bar{G}} \Delta \gamma  $ \\
\hspace*{5mm}$ \Xi_2 = \frac{2 \bar{G} \Delta \gamma}{||\bs{\tau}_{dev,n+1}^{trial}||} $  \\
\hspace*{5mm}$ \Xi_3 = \frac{1}{\Xi_0} - \Xi_2 + \Xi_1 $ \\
\hspace*{5mm}$ \Xi_4 = [\frac{1}{\Xi_0} - \Xi_3] \frac{||\bs{\tau}_{dev,n+1}^{trial}||}{\bar{G}} $\\

       Spatial elasticity tensor for hyperelasticity:\\
\hspace*{5mm}	$\bs{\mathbbm{c}}_{n+1}^{e,trial} = \frac{1}{J} [\bs{\mathbbm{c}}_{vol,n+1}^{e,trial} + \bs{\mathbbm{c}}_{dev,n+1}^{e,trial}]$ \\
\hspace*{5mm}          $\bs{\mathbbm{c}}_{dev,n+1}^{e,trial} = 2 \bar{G} ( \bs{I} - \frac{1}{3}\bs{I} \otimes \bs{I} )- \frac{2}{3}[(\bs{\tau}_{dev,n+1}^{trial} \otimes \bs{I}) + (\bs{I} \otimes \bs{\tau}_{dev,n+1}^{trial}  )]$ \\
\hspace*{5mm}           $\bs{\mathbbm{c}}_{vol,n+1}^{e,trial} = K [J^2 \bs{I} \otimes \bs{I} - (J^2-1) \bs{II}]$ \\

    Consistent tangent moduli:\\
\hspace*{5mm}$\bs{\mathbbm{c}}_{n+1}^{ep}= \frac{1}{J}[ \bs{\mathbbm{c}}_{vol,n+1}^{e,trial} + \bs{\mathbbm{c}}_{dev,n+1}^{e,trial}(1-\Xi_2) - 2 \bar{G} \Xi_3 \bs{n} \otimes \bs{n} - 2 \bar{G} \Xi_4 (\bs{n} \otimes \text{dev}[\bs{n}^2] )^{sym} $\\

     \bottomrule[2 pt]
    \end{tabular}
  \end{center}
\end{table}
\end{appendices}

\newpage
\bibliographystyle{unsrt}

\bibliography{References}

\end{document}